\documentclass{article}
\usepackage{a4wide}
\usepackage[pdftex]{graphicx}
\usepackage{amsfonts}
\usepackage{rotating}

\title{Numerical modelling of convection-diffusion-reaction problems with free boundary in 1D} \author{G.\ Ka\v{c}urov\'{a}\thanks{ Faculty of Mathematics and Physics. Comenius
University Bratislava, Mlynska dolina, 84215 Bratislava, Slovakia.
(E-mail {\tt gkacurova@yahoo.com}.)}}
\date{~}

\begin{document} \maketitle

\begin{abstract} 
We discuss a numerical method for convection-diffusion-reaction problems with a free boundary  in 1D. The method is based on the numerical modelling of the interface evolution, the transformation to a fixed domain problem and the approximation by an ODE
system. The interface evolution is modelled by means of the local shape of the corresponding travelling wave solution.
The method can be applied to many free boundary problems with a finite speed of the interface. The presented method can also approximate  some problems with an infinite speed of the interface for damped travelling wave type solutions. In the numerical experiments we compare our numerical solution with the analytical ones for some problems. 
\end{abstract}
\section{Introduction}
We shall consider convection-diffusion-reaction equations of the type
\begin{equation}
\label{1}
\partial_tu=\partial_x^2 u^n+b_0\partial_xu^{\gamma}+c_0u^m
\end{equation}
on the halfline $x\in (0,\infty)$ for  $t>0$, with the boundary conditions
\begin{equation}
\label{2}
a)\quad u=\phi(t),\  \mbox{or}\quad b)\quad -\partial_xu^n+b_0u^{\gamma}=q(t).
\end{equation}
at the boundary $x=0$. Here, $b_0$ and $ c_0$ are constants.
The initial condition is of the form
\begin{equation}
\label{3}
u(x,0)=u_0(x) \geq 0,
\end{equation}
where $u_0(x)> 0$ for $x\in (0,L_0)$ and $u_0(x)= 0$ for $x>L_0$. \
We are only looking for nonnegative solutions since negative solutions have no physical interpretation.
We assume that the parameters $n,\gamma$ and $m$ are nonnegative real numbers and guarantee the existence of a function $s(t)$ such that the solution is positive in $(0,s(t))$ and $u(x,t)=0$ for $x>s(t)$. The function $s(t)$ represents the time evolution of the interface. Such solutions have been intensively studied in the last 30-years in many papers, e.g., O.A. Olejnik, A.S. Kalashnikov, Czou-Ju-Lin \cite{[O-K-C]}, 
G.J.Barenblatt \cite{[B]}, D.G. Aronson \cite{[A]}, B.H.Gilding \cite{[G1],[G2]}, L.A.Peletier \cite{[P]}, R.Kersner \cite{[RK]}, B.H.Gilding and R.Kersner \cite{[G-K3],[G-K4]} etc. (see also the list of references in the previous papers).\\
We shall discuss the numerical approximation for (\ref{1})-(\ref{3}) based on interface modelling.
 In these
problems the free boundary $s(t)$ is not given explicitly in the model and its determination is included in the solution. Solutions to the linear convection-diffusion-adsorption problems ($n=m=\gamma=1$) in (\ref{1}) do not show the interface property. There, an initially localized support of the initial profile spreads with an infinite speed ($\dot{s}(t)=\infty$). Degeneration of the elliptic part ($n>1$)  is the main reason for the appearance of the interface.
However, also in the regular case ($n=1$)  the interface will appear under
special properties of the adsorption represented by $c_0$ and $m$. Very substantial results (necessary and sufficient conditions for $n,m,\gamma$ in (\ref{1})) to guarantee the appearance of the interface have been obtained in 
many papers, e.g. B.Song (\cite{[S]}), B.H.Gilding and R.Kersner (\cite{[G-K3],[G-K4]}), see also list of references there. If $b_0\neq 0, c_0<0$ and $ m,n,\gamma >0$  the existence of an interface 
is guaranteed when $\min(n,\gamma)> \min(m,1)$. This condition is also necessary (see \cite{[S]}).
If $b_0=0$ and $c_0=0$ in (\ref{1}) then the condition $n>1$  is the real criterion whether or not the interface appears. An equation with this property is called ``porous media equation" and an analytical solution has been found by Barenblatt \cite{[B]}. This equation is a mathematical model for adiabatic infiltration of a gas
into a porous media. Its solution is not a classical one (but a generalized, weak solution, defined by a corresponding integral identity) since $\partial_xu$ is not continuous along the interface, i.e. $x=s(t)$. The Barenblatt-Pattle solution of the Cauchy problem for the porous media equation
\begin{equation}
\label{5}
\partial_tu=\partial_x^2u^n, \quad n>1,\quad x\in (-\infty,+\infty)
\end{equation}
is of the form
\begin{equation}
\label{6}
u(x,t)=\left\{ 
\begin{array}[a]{c c}
\frac{1}{s(t)}[1-(\frac{x}{s(t)})^2]^{1/(n-1)},& \mbox{for}\ |x|<s(t)\\
0, & \mbox{for}\ |x|\geq s(t)
\end{array}
\right.
\end{equation}
where
$$ s(t)=[\frac{2n(n+1)}{n-1} (t+1)]^{1/(n+1)} $$
and the corresponding initial condition is
\[ u(x,0)=\left\{ 
\begin{array}[a]{c c}
\frac{1}{s(0)}[1-(\frac{x}{s(0)})^2]^{1/(n-1)},& \mbox{for}\ |x|<s(0)\\
0 ,& \mbox{for}\ |x|\geq s(0) .
\end{array}
\right.
\]
This analytical solution stimulated the theory of existence and uniqueness for the 
more general  ``porous media" type problems. More detailed information on this topic can be found in \cite{[G-K3],[G-K4]} and the references therein.\\
The solution of (\ref{5}) for $n=6$ and $ t \in (0,200)$ is drawn in Fig.1 (in $10$ equidistant time sections) and the evolution of its interface is drawn in Fig.2. We can readily verify the validity of the equation
\begin{equation}
\label{Bi}
\dot{s}(t)=-\frac{n}{n-1}\partial_xu^{n-1}|_{x\nearrow s(t)},
\end{equation}
which is the governing ODE model for this iterface. 
Solutions of the problem (1)-(3)  have the shape of damped travelling waves with eventually
sharp fronts at the interface. This makes the numerical approximation very difficult (the movement of the interface is not priori known). There exist numerical methods for free boundary problems that estimat the position of the interface and then apply a local adaptation of the space discretization, which are very expensive but have a limited precision (see e.g. \cite{[J-K],[K-M]} etc.). Some other results related to our method can be found in \cite{[BBChP],[Ch]}.\\

The main goal of this paper is the construction of an efficient and precise numerical method for solving  ``porous media" type problems (\ref{1})-(\ref{3}) . This method is based on the modelling (in terms of an ODE) of the time evolution of the interface $s(t)$ , which  significantly depends on the local profile of the solution near the interface. The unknown $s$ will be explicitly included into the solution of the problem. Then, we shall look for the solution in the interval $(0,s(t))$ only  (instead of $(0,L)$, where $L>L_0$). Consequently, we will transform (\ref{1}) to the fixed domain $(0,1)$ (using $y=\frac{x}{s(t)}$ ) and invoke a  special space discretization (finer in the neighbourhood of the
point $y=1$ which corresponds to the interface $x=s(t)$). This approximation of (\ref{1}) results in a corresponding system of ODEs including the  uknown $\dot{s}(t)$.
This system is completed with the ODE model for the interface  (of the same type as (\ref{Bi}) for (\ref{5})).  The final ODE system is solved with
an appropriate solver for stiff ODE systems.  Thus, the modelling of the interface evolution is the crucial point for the construction of our numerical method.  In this way we can approximate the sharp front of the solution at the interface  with higher accuracy. In our numerical experiments we have used 
the solver 'ode15s' from MATLAB library. \\
In Section 1 we describe the method in more detail. In Section 2 we construct the interface model for (\ref{1}) under some assumptions on parameters $n,m$ and $\gamma $.
In section 4 we construct the numerical approximation for (\ref{1}). In Section 5 we extend the results from the previous sections to the more general PDE
\begin{equation}
\label{4}
\partial_tu=\partial_x^2 a(u)+\partial_xb(u)+c(u)
\end{equation}
satisfying $a(0)=b(0)=c(0)=0$, $a'(0) = 0$ and
 preserving similar interface property of the solution as (\ref{1}) .  Also more general mathematical models involving $a(t,u),b(t,x,u),c(t,x,u)$ can be included too. These types of equations describe a large part of convection-diffusion-adsorption problems with many important practical applications.
In Section 6 we discuss the interface modelling in the non degenerate case $a'(0) > 0$.
These equations are of the ``oxygen-consumption" type, where appearence of the interface 
is generated by strong adsorption represented by $c(u)$.
In Section 7 we present some numerical experiments and we present comparisons with the analytical
solutions.  

\section{ Idea of the numerical approximation}

For a good numerical approximation of (\ref{1}) it would be desirable to use space discretization only in the moving domain  $(0,s(t))$ and to use also  moving
grid points which follow the interface and which are more dense at the sharp fronts
near the interface - see Fig.1. For this purpose we formally consider (\ref{1}) in the 
domain $(0,s(t))$ which by Landau's transformation, $y=\frac{x}{s(t)}$, is transformed it to the
fixed domain $ (0,1)$. Then, for $\bar{u}(y,t):=u(x,t)$ we obtain the equation
\begin{equation}
\label{tr}
 \bar{u}_t=\bar{{\cal L}} (\bar{u})- \frac{\dot{s}}{s} y \partial_y\bar{u},\ y\in (0,1)
\end{equation}
where ${\cal L}(u)$ is the rhs of (\ref{1}) and $\bar{{\cal L}} (\bar{u})$ is its transform in the $y$-variable. If we can model $\dot{s}$ in terms of $\partial_y\bar{u}(1,t)$, resp. $\partial_y^2\bar{u}(1,t)$ (which describe local properties of $u$ in a 'small' neighbourhood of $s(t)$), e.g., in the form (see (\ref{1}))
\begin{equation}
\label{int}
 \dot{s}= {\cal B}(\partial_y\bar{u}(1,t)),\ \mbox{resp.}\    \dot{s}={\cal B}(\partial_y\bar{u}(1,t),\partial_y^2\bar{u}(1,t)), 
\end{equation}
then we can eliminate $ \dot{s}$ in (\ref{tr}) by means of the rhs of (\ref{int}).
Then the complete system consists from (\ref{tr}) and (\ref{int}). We approximate it
using the space discretization $\{y_i\}_{i=1}^N \in (0,1)$ of the $y$-variable. The resulting system of ODEs is stiff and we solve it by a corresponding ODE solver. Now, the unknown interface is explicitly included in the solution. The grid points
$\{y_i\}_{i=1}^N $ correspond to the moving grid points $ \{x_i(t)\}_{i=1}^N \in (0,s(t))$, with $x_N(t)=s(t)$ for (\ref{1}). 

The method used cannot be  extended to more space dimensions  (except of radial symmetric ones).
The mathematical model for the interface movement (its normal velocity) depends not only on local properties
of the solution (e.g., the normal space derivative) but also on the geometry of the support of the solution, which in 1D is very simple.

\section{ Modelling of the interface}

 The appearance of the interface in (\ref{1}) is closely related to the existence of the travelling wave solution (see \cite{[G-K3],[G-K4]} ) of the form
$u(x,t)=f(\xi),\ \xi=x-\sigma t$ where 
\[ \left\{
\begin{array}[a]{c c}
f(\xi)=0 ,& \mbox{for}\ \xi\ge \xi_0\\
f(\xi)>0, & \mbox{for}\ \xi < \xi_0	
\end{array}
\right.
\]
for some $\xi_0\in (-\infty,+\infty)$. The function $f$ has to be determined from the ODE
$$
-\sigma f'=nf^{n-1}f''+n(n-1)f^{n-2}(f')^2+b_0\gamma f^{\gamma-1}(f))'+c_0f,\quad \xi \in (-\infty,\infty).
$$
The function $\theta$ defined as 
$$\theta(f(\xi))=-(f^n)'(\xi), $$
has to satisfy the nonlinear Volterra equation (see \cite{[G-K3]})
\begin{equation}
\label{travel}
 \theta(s)=\sigma s +b_0s^{\gamma} -nc_0\int_0^s \frac{r^{m+n-1}}{\theta(r)}dr.
 \end{equation}
In  \cite{[G-K3]} (Lemma 13) one proves: There exists a solution $\theta$ of  (\ref{travel}) for all  $\sigma >\sigma^{\ast}$ (for suitable  $\sigma^{\ast}>0$) provided 
the conditions (i)-(iii) below with parameters $n\ge 1,\  m>-n$  and $ \gamma \ge 1$  are fulfilled. Moreover, there exists a $\sigma^{\ast \ast}> \sigma^{\ast}$ so that
 \begin{equation}
\label{ftw}
 \theta(s)\approx \theta_0s^{q}\quad \mbox{as}\  s\searrow 0
\end{equation} 
for some $\theta_0$ and $ q>0$ depending on $n,m,\gamma$. These conditions are
\[
\begin{array}[a]{c c c c}
(i)& c_0<0 ;&\quad ( \mbox{in this case}\ &  q=\min \{(n+m)/2,1\})\\
(ii)& c_0=0 ;&\quad (\mbox{in this case}\ &  q=1)\\
(iii)& c_0>0,\ m+n\ge 2;&\quad ( \mbox{in this case}\ &   q=1).
\end{array}
\]
The conditions (i)-(iii) are also necessary for the existence of the solution of (\ref{travel}).\\
The case $0<n<1$ is also discussed in \cite{[G-K3]} (Lemma 13), but we cannot aply it in our analysis.\\
The relation (\ref{ftw}) is equivalent to
$$ \theta(s)= \theta_0s^{q}+o(s^q),\quad \mbox{with}\ s^{-q}o(s^q)\to 0,\ \mbox{for}\ s\searrow 0.$$
As a consequence of (\ref{ftw}),
also the front of the semi-wave $f$  (localized at 
the point $\xi_0$), satisfies the approximation (see Remark 1 below)
\begin{equation}
\label{ftw2}
f(\xi)\approx d^{\beta} (\xi_0-\xi)^{\beta}\quad  \mbox{for}\ \xi\nearrow \xi_0 .
\end{equation} 
We shall use this property in the modelling of $\dot{s}(t)$ in terms of local properties of the solution
at the interface (e.g. $\partial_xu(s(t),t)$) and model parameters. This will be substantially used in our numerical method. In fact, we shall assume the asymptotic property (\ref{ftw2}) to hold also for the first and second derivatives of $f$ (see (\ref{ftw3}) below).\\
{\bf Remark 1}. From (\ref{ftw})and $\theta(f(\xi)):=\frac{d}{d\xi}f^n(\xi))$ (provided $f$ is smooth, decreasing and $f(\xi)=0$ for $\xi\ge {\xi}_0$ - ``travelling semi-wave") it follows 
$$f'(\xi)\approx \frac{\theta_0}{n}f(\xi)^{\alpha},\quad  \alpha=q-n+1$$
and consequently, by integration  over $(\xi,{\xi}_0)$, we obtain (see (\ref{ftw2})),
 $$ f(\xi)\approx [\frac{\theta_0(n-q)}{n}]^{\beta}{(\xi}_0-\xi)^{\beta},\quad \mbox{where}\ \beta=\frac{1}{n-q}, d=\frac{\theta_0(n-q)}{n}$$
 under the assumption $n\ne q$.\\

 The development of an interface model of the type (\ref{int}) in a rigorous mathematical way requires a very deep analysis. Moreover, there are no  results available in this respect (except of a few, very special cases). We will construct such a model under very strong smoothness assumptions on the solution in the interval $(0,s(t)]$ (up to the free boundary) and the assumption (\ref{ftw2}) concerning the shape of travelling semiwave corresponding to (1). In fact we assume (\ref{ftw2})
in its stronger form
$$ f(\xi)\approx d^{\beta} (\xi_0-\xi)^{\beta},\ f'(\xi)\approx \beta d^{\beta} (\xi_0-\xi)^{\beta-1}  \quad  \mbox{for}\ \xi\nearrow \xi_0 ,$$
\begin{equation}
\label{ftw3}
 [f^n]'(\xi)\approx d^{n\beta} n\beta (\xi_0-\xi)^{n\beta-1}\quad  \mbox{for}\ \xi\nearrow \xi_0,
\end{equation}
$$ [f^n]''(\xi)\approx d^{n\beta} n\beta (n\beta-1)(\xi_0-\xi)^{n\beta-2}\quad  \mbox{for}\ \xi\nearrow \xi_0$$
where $\beta$ and $\ d$\ are unknown parameters to be determined.\\

Our main result in this section is\\
{\bf Theorem 1}. Suppose (\ref{ftw3}) and $n> 1,\ 0\le m <1,\gamma \ge 1$ . When the solution of (1) is regular up to the interface (the equation is satisfied in the limit sense at the interface), then the interface $s$ of (1) is governed by ODE
\begin{equation} 
\label{14}
\dot{s}(t)=Q(w)=\left\{
\begin{array}[a]{c c}
- \frac{n}{(n-1)}\partial_xw +\frac{c_0}{\partial_xw}-G(\gamma)b_0&\  \mbox{for}\ m+n=2;\\ 

- \frac{n}{(n-1)}\partial_xw -G(\gamma)b_0 &\  \mbox{for}\ m+n>2\\
- \frac{n}{(n-1)}\partial_xw -G(\gamma)b_0 &\  \mbox{if}\ c_0= 0

\end{array}
\ ,\quad   x=s(t) 
\right.
\end{equation}
where $\beta=\frac{1}{n-1}$,\  $u=w^{\beta}$  and
$G(r)=0$ for $r>1$,\ $G(1)=1$. \\
{\bf Proof}.
 Our starting point is the condition $u(s(t),t)=0$ for all $t\in (0,T)$. Differentiating this equation with respect to $t$ we obtain formally
\begin{equation}
\label{12}
 \dot{s}(t)=- \frac{\partial_tu}{\partial_xu}\ \mbox{in the sense}\ x\nearrow s(t),
 \end{equation}
which cannot be used practically, since usually we have $\partial_xu \to -\infty$ for $x\nearrow       s(t)$. But invoking the assumption (\ref{ftw3}),  we can conclude that $\partial_xu^{\frac{1}{\beta}}$ is finite for $x\nearrow s(t)$ for some $\beta >0$ . Using the smoothness of the solution (the equation (1) holds up to the interface $x=s(t)$) we can eliminate $\partial_tu$ in  (\ref{12}) by the rhs of (1) and then we insert $u=d^{\beta}(s(t)-x)^{\beta}+o_0$
  in the small neighbourhood of the interface (for $x\nearrow s(t)$). Due to (\ref{ftw3}), we obtain
\begin{equation}
\label{13}
\dot{s}(t)=\lim_{x\nearrow s}\frac{1}{d^{\beta}(s-x)^{\beta-1}+o_1}\{ d^{n\beta}n\beta( n\beta-1)(s(t)-x)^{n\beta-2}-
\end{equation} 
$$b_0d^{\gamma \beta}\gamma \beta (s-x)^{\gamma \beta-1} +c_0d^{m\beta}(s-x)^{m\beta}+ o_2+ o_3+ o_4 \}$$ 
with the lower order terms represented by $o_1-o_4$ ( because of (\ref{ftw3})), where
 $$ o_0\equiv o((s-x)^{\beta}),\ o_1\equiv o((s-x)^{ \beta-1}),\  o_2\equiv o((s-x)^{n\beta-2}), $$   $$\ o_3\equiv o((s-x)^{\gamma \beta-1}), \ o_4\equiv o((s-x)^{m \beta}).$$ 
The assumtions on $n,m$ and $\gamma$ guarantee the existence of the interface, i.e. $|\dot{s}(t)|<\infty$ because of  \cite{[G-K3]} (Lemma 13). Our goal is to find such a $\beta>0$  that this will be achieved.
There are more combinations of parameters $n,m,\gamma$ to guarantee that such a $\beta$ could be found. Deviding the nominator  on the  rhs in (\ref{13}) by $(s-x)^{ \beta-1}$, we obtain three terms containing $s-x$ with the exponents $(n-1)\beta -1;\ (\gamma-1)\beta ;\ (m-1)\beta +1$ (corresponing to diffusion, convection and reactive part) . To guarantee $|\dot{s}(t)|<\infty$, we have to choose $\beta >0 $ such that all three exponents are nonnegative and at least one of them has to be zero.
The case $\beta \le 0 $ doesn't correspond to the required property of the semiwave. Then, for $n+m\ge 2$ the choice $\beta =\frac{1}{n-1} $ fulfilles our requirements.
This choice also includes the case $c_0=0$ .  
 Then from 
(\ref{13}) we deduce the  formula for $\dot{s}(t)$ which contains the uknown
 parameter $d$ -see (\ref{ftw2}).
In the following we shall determine  $d$  . For this purpose consider the new unknown $w$ defined by
$u=w^{\beta}$  and rewrite (1) in terms of $w$. In the small neighbourhood of the interface  ($x\nearrow s(t)$) we have
$w(x,t)\approx d(s(t)-x)$ due to (\ref{ftw3}) and hence 
$$|\frac{w(x,t)-w(s(t),t)}{s(t)-x}-d| \to 0,\ \mbox{as} \ x\nearrow s(t).$$
 From this we conclude that $\partial_xw\to -d$ when $x\nearrow s(t)$.\\
 Now, we substitute $d=-\partial_xw$ into (\ref{13}) and we obtain the required ODE model for the interface.\\

{\bf Remark 2}. It is very difficult to verify the assumptions (\ref{ftw3}). However
 consequences from them,  leading to the interface models, are verified  in some examples in Section 7, where the analytical solution is available. Since assumptions (\ref{ftw3}) could be violated in some model problems, our modelling of the interface should be justified numerically in practical applications. \\
{\bf Remark 3}. The analysis used in choosing $\beta$ so that $|\dot{s}(t)|<\infty$  is rough and cannot distinguish the sign of $b_0$ and $c_0$ which is important for the existence of the interface. The analysis for the adsorption and reaction is the same. From (i) and (iii) it follows that the interface exists also in the case $n+m<2$ provided $c_0<0$ (i.e. in the adsorption case), while for $c_0>0$ the interface doesn't exist. \\
{\bf Remark 4}.  In the case $n+m \ge 2$ it follows from (i),(iii) that $q=1$. Then, our $\beta$ obtained in the proof of Theorem 1 coincides with the one from Remark 1 what supports our arguments. In fact, the assumptions in Theorem 1 include also the additional solution $\beta=\frac{1}{1-m}$ in the case $n+m>2$. However, this $\beta $ is excluded by the fact $q=1$ (see (i), (iii) and Remark 1). It means that our analysis leading to Theorem 1 is not satisfactory for the unique determination of $\beta$.\\

\section{Numerical approximation  of (1)}
The mathematical model for the time evolution of $s$ given by Theorem 1 is of practical use, since $\partial_xw$ is finite and nonzero. From this reason we also rewrite (1) in  terms of $w$ using the transformation $u=w^{\beta}$ where $\beta =\frac{1}{n-1}$ . Then, we obtain
\begin{equation}
\label{15}
\partial_tw=n(n\beta-1)w^{(n-1)\beta}\partial_x^2w+nw^{(n-1)\beta-1}(\partial_xw)^2+
b_0\gamma w^{(\gamma-1)\beta}\partial_xw+\frac{c_0}{\beta} w^{(m-1)\beta+1}.
\end{equation}
 We have to solve this equation in the moving domain $(0,s(t))$. We transform (\ref{15}) to the fixed domain $(0,1)$ using  Landau's transformation $ y=\frac{x}{s(t)}$, with $\bar{w}(y,t)=w(x,t)$.
 For the sake of simplicity, in the following we write $w$ in the place of $\bar{w}$. Then, we obtain
$$
\partial_tw=n(n\beta-1)\frac{w^{(n-1)\beta}}{s^2}\partial_y^2w+
n\frac{w^{(n-1)\beta-1}}{s^2}(\partial_yw)^2+
b_0\gamma \frac{w^{(\gamma-1)\beta}}{s}\partial_yw+$$
 $$+\frac{c_0}{\beta} w^{(m-1)\beta+1}+ y\frac{\dot{s}}{s}\partial_yw.$$
This equation has to be completed with the ODE for $\dot{s}$ given by Theorem 1 and $\dot{s}$ at the rhs has to be eliminated by means of the rhs of (\ref{14}) . Then, we have to solve the coupled PDE and ODE system 
 \begin{equation}
\label{16}
\partial_tw=n(n\beta-1)\frac{w^{(n-1)\beta}}{s^2}\partial_y^2w+
n\frac{w^{(n-1)\beta-1}}{s^2}(\partial_yw)^2+
b_0\gamma \frac{w^{(\gamma-1)\beta}}{s}\partial_yw+$$
 $$+\frac{c_0}{\beta} w^{(m-1)\beta+1}+ y\frac{Q(w)}{s}\partial_yw,
\end{equation}
\begin{equation}
\label{17}
\dot{s}=Q(w)
\end{equation}
on a fixed domain $y\in (0,1),\ t>0$, where $Q$ is from (\ref{14}) (there, $\partial_xw$ has to be replaced by $\frac{\partial_yw}{s}$). 
The PDE  (\ref{16}) at the point $(y,t)$ depends also on the point $(1,t)$ because of $Q(w)$. Consequently, it is non-local.
 In the numerical approximation of (\ref{16})-(\ref{17})
we use  space discretization (method of lines) which results in an ODE system.
We  generally use a nonequidistant partition of $(0,1)$. Let $\alpha_i>0,\ \forall i=1,...,N$,\  $\alpha_0=0$ and consider the grid points $y_i=\sum_{j=0}^i \alpha_j,\ \forall i=0,...,N$. We approximate $w(y_i,t)\approx C_i(t)$. Denote by $\ell_i(y)$ the Lagrange polynomial of the second order through the points $(y_{i-1},C_{i-1})$,\ $(y_{i},C_{i})$ and  $(y_{i+1},C_{i+1})$. We approximate $\partial_y^2w(y_i,t)$ and  $\partial_yw(y_i,t)$ by the corresponding expressions given by
$\frac{d^2 \ell_i}{dy^2}|_{y=y_i}$ and $\frac{d \ell_i}{dy}|_{y=y_i}$.\\
We approximate the derivative
$\partial_yw(1,t)$ appearing in $Q$  by
$$\partial_yw(1,t)=\frac{d \ell_B}{dy}|_{y=1},$$
where $\ell_B$ is the second order Lagrange polynomial throudh the points $(y_{N-2},C_{N-2})$ and   $(y_{N-1},C_{N-1})$,\ $(1,0)$.
Our ODE system now reads
 \begin{equation}
\label{18}
\dot{C}_i=n(n\beta-1)\frac{C_i^{(n-1)\beta}}{s^2}\frac{d^2 \ell_i}{dy^2}|_{y=y_i}+n\frac{C_i^{(n-1)\beta-1}}{s^2}(\frac{d \ell_i}{dy}|_{y=y_i})^2+b_0\gamma \frac{C_i^{(\gamma-1)\beta}}{s}\frac{d \ell_i}{dy}|_{y=y_i}+
\end{equation}
 $$+\frac{c_0}{\beta} C_i^{(m-1)\beta+1}+
y_i\frac{Q(C)}{s}\frac{d \ell_i}{dy}|_{y=y_i},$$
 \begin{equation}
\label{19}
\dot{s}=Q(C),\quad \mbox{where}\ \partial_yw(1,t)\leftrightarrow  \frac{d \ell_B}{dy}|_{y=1}
\end{equation}
  for $i=1,...,N-1$, with $\beta=\frac{1}{n-1}$. In the case of Dirichlet boundary conditions this is the resulting system of ODE,
where we use $C_0(t)\equiv \phi(t)$ (see (\ref{2}),\ case a)). In the case of a Robin condition (see (\ref{2}),\ case b)), the system (\ref{18})-(\ref{19}) has to by completed by the additional ODE  at the point $y=0$. In this case we make use of a ghost point $y_{-1}=-\alpha_1$. The missing value $C_{-1}$ is determined
from the boundary condition $(\ref{2})$, (case b), by means of $q$ and $C_1$ from the relation
$$ C_{-1}=C_1+ s(-b_0C_0^{\gamma \beta}+q)\frac{2\alpha_1}{n\beta C_0^{n\beta-1}}.$$
Then, the additional ODE in $y=0$ is the same as (\ref{18}) for $i=0$ using  $C_{-1}$.\\
In some physical problems it is important to know the support of the solution
 which in 1D is given by the position of the interface. This is explicitly included in our approximation represented by (\ref{18})- (\ref{19}). This gives us the possibility to determine $s(t)$ with higher accuracy then the approximations without interface modelling.\\
{\bf Remark 5}.
The system (\ref{18})-(\ref{19}) is stiff and a corresponding ODE solver has to be used. This system is a good starting point to solve the original free boundary value problem. The first advantage is
that we approximate the solution on its support only. The second one lies in the  space discretization. In fact, we have moving grid points (in the original\ $x$- variable) by means of which sharp fronts near the interface can be approximated very well for all $t$. We shall demonstrate this
in our numerical experiments presented in  Section 7 . The presented numerical approximation requires positive initial profile on its support $(0,s(0))$. If our initial profile is zero (provided nonzero boundary condition at $x=0$ is considered), then we have to start with the ``small'' positive auxiliary initial profile on the
`` small ''  domain $(0,s(0))$.
 Also some compatibility between the initial and boundary conditions (at the point $(t,x) =(0,s(0))$) help to avoid numerical instabilities in the ODE solver at the starting time point. This ``small'' regularization does not change the expected solution significantly. The efficiency of the numerical method considerably depends on the strategy of grid points  distribution. This allows reducing the size of the system (\ref{18}) without decreasing the accuracy.\\
The used space discretization of the  PDF and its approximation it by an  ODE- system is of the second order and the order of ODE solver is regulated by its tolerance parameters. 

\section {Extension of the method to PDE (\ref{4})}

From the theory of porous media type solutions one can notice that the evolution of the front of the solution substantially depends on the order of degeneration and adsorption and on local properties of the solution. Therefore, the same idea of the numerical approximation used on (\ref{1}) can  also  be applied in a more general case of (\ref{4}), where we replace $a(s) \leftrightarrow a(t,s),\ b(s)\leftrightarrow b(x,t,s),\ c(s)\leftrightarrow c(x,t,s)$, provided $a,b$ and $c$  are asymptotically polynomial with respect to $s$ 
at the point $s=0$. More precisely, we assume
\begin{equation}
\label{20}
a(t,u)=g(t,u)u^n,\ b(x,t,u)=b_0(x,t)u^{\gamma},\ c(x,t,u)=c_0(x,t)p(u)u^m  
\end{equation}
$$\mbox{with}\quad g(t,0)\neq 0,\ p(0)\neq 0\quad  \mbox{and}\ g\in C^2.$$
Modelling the interface in this setting (generalization of (\ref{4})), we proceed similarly as in Section 2 and obtain.\\
{\bf Theorem 2}. Let the assumptions of Theorem 1 hold true. Moreover, suppose (\ref{20}).
Then, the governing ODE for the interface in the  genaralized case of (\ref{4}) is given by
\begin{equation}
\label{21}
\dot{s}(t)=Q(w)=\left\{
\begin{array}[a]{c c}
- \frac{n g(t,0)}{(n-1)}\partial_xw +\frac{c_0(s(t),t)p(0)}{\partial_xw}-G(\gamma)b_0(s(t),t)&,\  \mbox{for}\ m+n=2;\\ 
- \frac{n}{(n-1)}\partial_xw -G(\gamma)b_0(s(t),t) &,\  \mbox{for}\ m+n>2,\\
- \frac{n}{(n-1)}\partial_xw -G(\gamma)b_0(s(t),t) &,\  \mbox{for}\ c_0=0,
\end{array}
\right.
\end{equation}
where  $u=w^{\beta}$ with $\beta=\frac{1}{n-1} $ and $G(r)=0$ for $r>1$,\ $G(1)=1$ and $u=w^{\beta}$ with $\beta=\frac{1}{n-1} $.\\
Using this ODE model for the interface, we approximate the considered PDE similarly as 
in the case of (\ref{1}).

\section{Modelling of the interface in the nondegenerated case  $\partial_u a(t,0)>0$.}
Cosider the generalization of (\ref{4}) as in the previous section with the nondegenerate elliptic part  $\partial_u a(t,0)>0$.
As it was mentioned in the introduction, the interface can appear also in the  nondegenerated case when adsorption is sufficiently strong for small values of $u$ . As an example we can take oxygen-consumption type  problem
\begin{equation}
\label{22}
\partial_tu=D\partial_x^2u-H(u),\ x\in (0,\infty),t>0
\end{equation}
with $u(x,0)=u_0(x)>0$ for $x\in (0,L_0)$, $u_0(x)=0$ for $|x|>L_0$, where $D$ is the diffusion coefficient and
$$ H(z)=\left\{
\begin{array}[a]{c c}
1 & \ \mbox{for} \ z>0 \\
\\
0 & \ \mbox{for}\ z=0
\end{array}
\right.
$$
represents the consumption of the oxygen with concentration $u$. This model problem has also an interface with finite speed and has been analysed some decades ago (see, e.g., \cite{[C]}).
There exists an analytical solution wich is a polynomial of the second order at the front (near the interface). 
This motivates us to look for a polynomial solution at the interface also in the more general case (\ref{4}) considered in the previous section. The order of this polynomial can be determined from the information that the considered problem has a finite speed of interface propagation. Then, for the interface modelling we proceed similarly as in the degenerated case. 
The existence of the interface (in case of (\ref{4})) has been established in
\cite{[G-K3]} (Theorem 12, Assertion (i)).\\
For simplicity we consider a generalization of the mathematical model (\ref{4}) (of oxygen-diffusion type)
under the assumptions
\begin{equation}
\label{23}
a=a_0(x,t)u,\ b=b_0(x,t)u,\ c=c_0(x,t)H(u)u^m  
\end{equation}
$$\mbox{with}\quad a_0(x,t)\ge \delta>0,\ c_0(x,t)\le -\delta,\ 0\le m<1.$$
A more special case when $c=c_0H(u)$ has been treated  numerically in \cite{[CK]} and we have extend this idea to the more general case of (\ref{4}), where (\ref{23}) is fulfilled. Another approximation method based on the solution of variational inequalities has been used in \cite{[B]}.\\
We start with the formal ``interface condition" $\dot{s}=-\frac{\partial_t u}{\partial_x u},$
where we substitute the rhs from (\ref{4}) with $\partial_t u$ under the smoothness assumption on the solution and with the validity of (\ref{4}) up to the interface. Then, we get
\begin{equation}
\label{24}
\dot{s}(t)=-\frac{\partial_x^2a(x,t,u)+\partial_xb(x,t,u)+c(x,t,u)}{\partial_xu},
\end{equation}
where the rhs is understood in the limit sense $x\nearrow s(t)$. Again, we consider the transformation $u=w^{\alpha}$ for some $\alpha>0$ in (\ref{4}) in the setting (\ref{23}). We obtain
\begin{equation}
\label{25}
\partial_t w =a_0[\partial_x^2w+(\alpha-1)\frac{\partial_xw)^2}{w}] +[b_0+2\partial_xa_0]\partial_xw+\frac{1}{\alpha} \frac{w^{(m-1)\alpha+2}c_0}{w}+\frac{1}{\alpha}(\partial_x^2a_0+\partial_xb_0)w
\end{equation}
where the variables $x$ and $t$ in functions $a_0,b_0$ and $c_0$ are omitted. We are looking for such an  $\alpha$  that (\ref{25}) holds  up to the interface $x=s(t)$.
Since $w(s(t),t)=0$ ($w(x,t)\to 0$ for $x\to s(t)$), it must hold that
\begin{equation}
\label{26}
(\alpha-1)a_0(\partial_xw)^2+\frac{c_0}{\alpha}w^{(m-1)\alpha+2} \to 0\quad \mbox{for}\ x\to s(t).
\end{equation}
Due to $a_0(s(t),t)\ge \delta >0$, we conclude
\begin{equation}
\label{27}
\partial_xw=- \sqrt{ -\frac{c_0}{\alpha(\alpha-1)a_0}w^{(m-1)\alpha+2}} \quad \mbox{for}\ x\nearrow s(t).
\end{equation}
The negative sign  is chosen due to the fact that $w$ is decreasing as $x\nearrow s(t)$.
Since 
\begin{equation}
\label{interface}
\dot{s}(t)=-\frac{\partial_tu}{\partial_xu}=-\frac{\partial_tw}{\partial_xw}
\end{equation}
and  $| \dot{s}(t)|<\infty$, we choose $\alpha $ so that $\partial_xw\neq 0$. Then, with respect to (\ref{27})
we choose  $$ (m-1)\alpha+2=0\quad \Rightarrow \ \alpha=\frac{2}{1-m}>1.$$
This guarantees that $|\partial_xw|<\infty$. To make use of (\ref{25}) in (\ref{interface}) we have to compute the limit
$$
Z_w=\lim_{x\nearrow s(t)}\left(\frac{(\alpha-1)a_0(\partial_xw)^2+\frac{1}{\alpha}c_0}{w}\right).
$$
Both, the nominator and denominator tend to  zero as $x\nearrow s(t)$, so the L'Hospital rule can be used. We obtain
\begin{equation}
\label{28}
Z_w=(\alpha-1)\partial_xa_0(s(t),t)\partial_xw + 2(\alpha-1)a_0(s(t),t)\partial_x^2w
+\frac{1}{\alpha} \frac{\partial_xc_0(s(t),t)}{\partial_xw}=
\end{equation}
$$ 2(\alpha-1)a_0(s(t),t)\partial_x^2w-$$ $$ (\alpha-1)\partial_xa_0(s(t),t)\sqrt{-\frac{c_0(s(t),t)}{\alpha(\alpha-1)a_0(s(t),t)}}-
\frac{1}{\alpha}\partial_xc_0(s(t),t)\sqrt{-\frac{\alpha(\alpha-1)a_0(s(t),t)}{c_0(s(t),t)}}.$$
Finally, from (\ref{24})-(\ref{28}) we conclude
\begin{equation}
\label{29}
\dot{s}(t)=(2\alpha-1)a_0(s(t),t)\sqrt{-\frac{\alpha(\alpha-1)a_0(s(t),t)}{c_0(s(t),t)}}\partial_x^2w(s(t),t)-b_0(s(t),t)+
\end{equation}
$$2\partial_x a_0(s(t),t)+(\alpha-1)\partial_xc_0(s(t),t)\frac{a_0(s(t),t)}{c_0(s(t),t)},\ \mbox{where}\
\alpha=\frac{2}{1-m}\ \mbox{and}\ u=w^{\alpha}. $$
{\bf Theorem 3}. Consider (\ref{4}) under the assumption (\ref{23}). Let the solution $u$ be $C^2$- smooth up to the interface $s(t)$. Then, the solution admits an interface $s(t)$ and this is governed by the ODE (\ref{29}).\\
In the case $m=0$, the solution $u$ is asymptotically a quadratic polynomial at the interface, since $u=w^2$, $\partial_xw\ne 0$, and $u$ is finite on the interface. This is the case for a classical oxygen-consumption problem (\ref{22}).\\
  This result can be extended to the more general case (\ref{4}), with 
  $a=a_0(x,t)g(u)u, b=b_0 u^{\gamma}$, where $g(0)>\delta >0,\ \gamma \ge 1$ and $c$ is as in (\ref{23}).\\
 The numerical approximation of the generalized oxygen-consumption problem (\ref{4}) under the
 conditions (\ref{23}) proceeds along the same lines as in Section 4, using our governing equation (\ref{29}). The approximation $\partial_x^2w$ on the interface by means of $\partial_y^2\ell_B|_{y=1}$ is generaly not satisfactory.  The third order polynomial through the points $$(y_{N-3},C_{N-3}),(y_{N-2},C_{N-2}),(y_{N-1},C_{N-1}),(1,0).$$
has to be used.

 \section{ Numerical experiments}

 To make an advantage of our numerical approximation a suitable discretization (in the $y\in (0,1)$-variable) can be chosen. Since $y=1$ corresponds to the interface, where the front of the travelling and damping semiwave is situated, we increase the density of
  grid points in the neighbourhood of  $y=1$. At the same time we decrease the density of grid points at the trailing part of the damping semiwave type solution to keep the  set of all grid points small. This of course depends on the solution profile and on the type of 
 the solved model. The distribution of the grid points plays an important role in the accuracy and efficiency of the numerical realization. Numerical approximation without the interface information requires much more grid points in the space discretization then our method with interface modelling to obtain the same order of accuracy. 
 For large $t$
 the profile of the solution can be significantly changed comparing with that one for small $t$. In that case we can restart the numerical solution of the corresponding ODE-system for large $t$ with a new grid (with a changed strategy of grid points distribution) and the corresponding initial profile.
 The main goal is to obtain the same required accuracy with a minimal set of grid points, which determines the dimension of our ODE system to be solved. We have used the following discretizations. The interval
 $(0,1)$ is divided equidistantly in $M$ subintervals and the last $d$ intervals are successively subdivided :to $2^p,3^p,...(d-1)^p$ additional subintervals for $p=1,2,3$. We denote these discretizations by $D_p$. Mostly, we will apply a smooth geometrical
 dicretization, where the length of the $(i+1)$-th subinterval is a $q<1$ factor smaller than that of the $i$-th.
 At this strategy, $q$ is determined from the prescribed number of all grid points $N$ and the length of the first subinterval $1/M$. We denote this discretization by $D_4$. Space discretization $D_4$ seems to be most efficient in most of our numerical applications.\\
 In our numerical experiments we solve the resulting ODE-system by means of solver ode15s from MATLAB library, which is developed for stiff ODE-systems.

\subsection{Barenblatt solution.} 

The interface model given by Theorem 1 (the case with $c_0=0$) coincides with that one in ({\ref{Bi}). 
In Fig. 1 we present our numerical solution using $n=6$ and the discretization $D_4$ with parameters $M=10$ and $ N=20$. The corresponding analytical solution cannot be distinguished graphically from the numerical one. The time evolution of the interfaces is shown in Fig. 2. The interfaces cannot be distinguished, too, even during the long time period $(0,200)$.\\
\begin{figure}
\begin{center}
  \setlength{\unitlength}{1cm}
\begin{picture}(12,8)
\put(0,0){\mbox{
    \includegraphics[width=12cm,height=6cm]{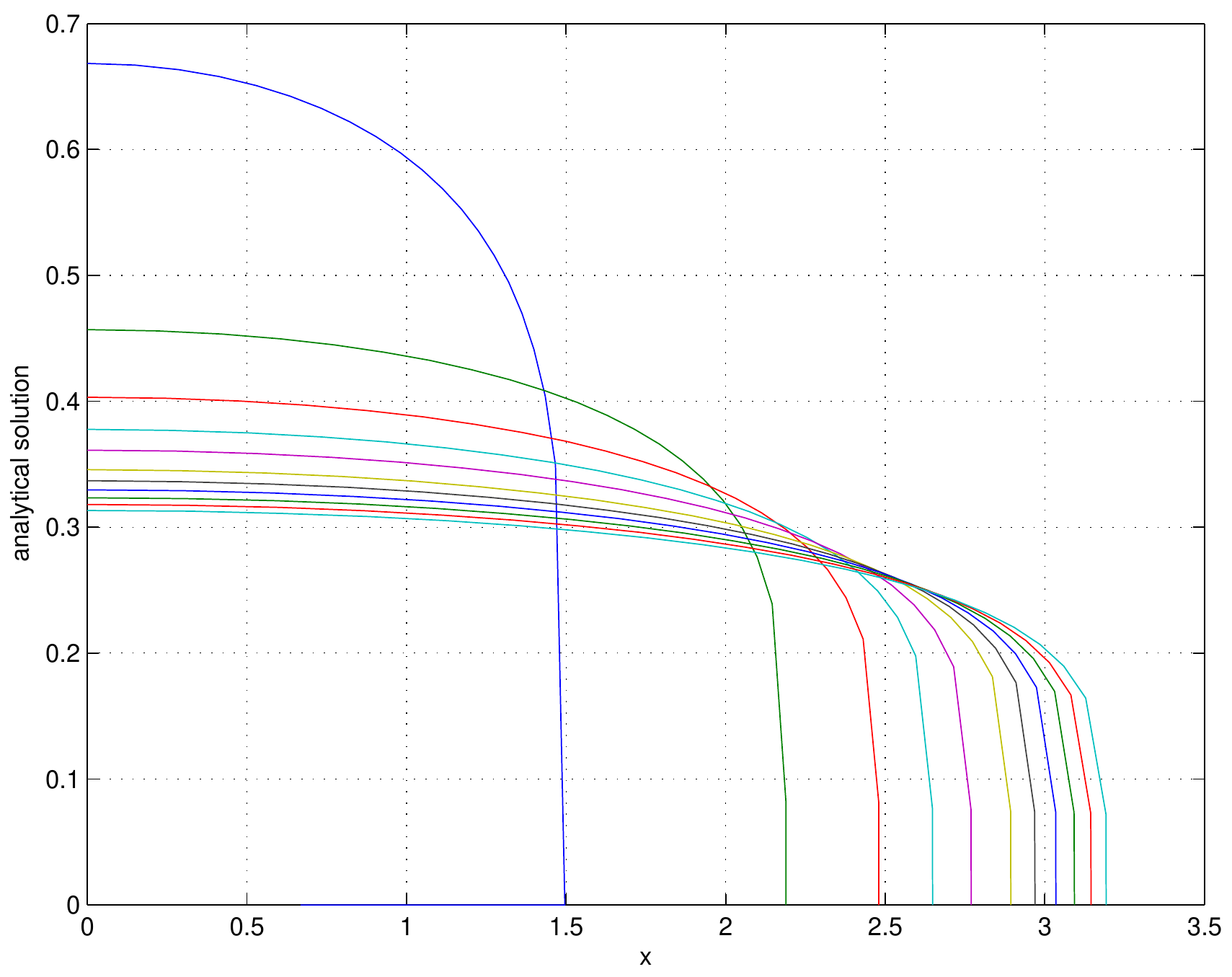}
    }}
\end{picture}
\caption{Numerical solution of (\ref{5}) with $p=6$ in $10$ equidistant time sections $t\in [0,200]$}   
\label{fig1}
\end{center}
\end{figure}

\begin{figure}
\begin{center}
  \setlength{\unitlength}{1cm}
\begin{picture}(12,8)
\put(0,0){\mbox{
    \includegraphics[width=12cm,height=6cm]{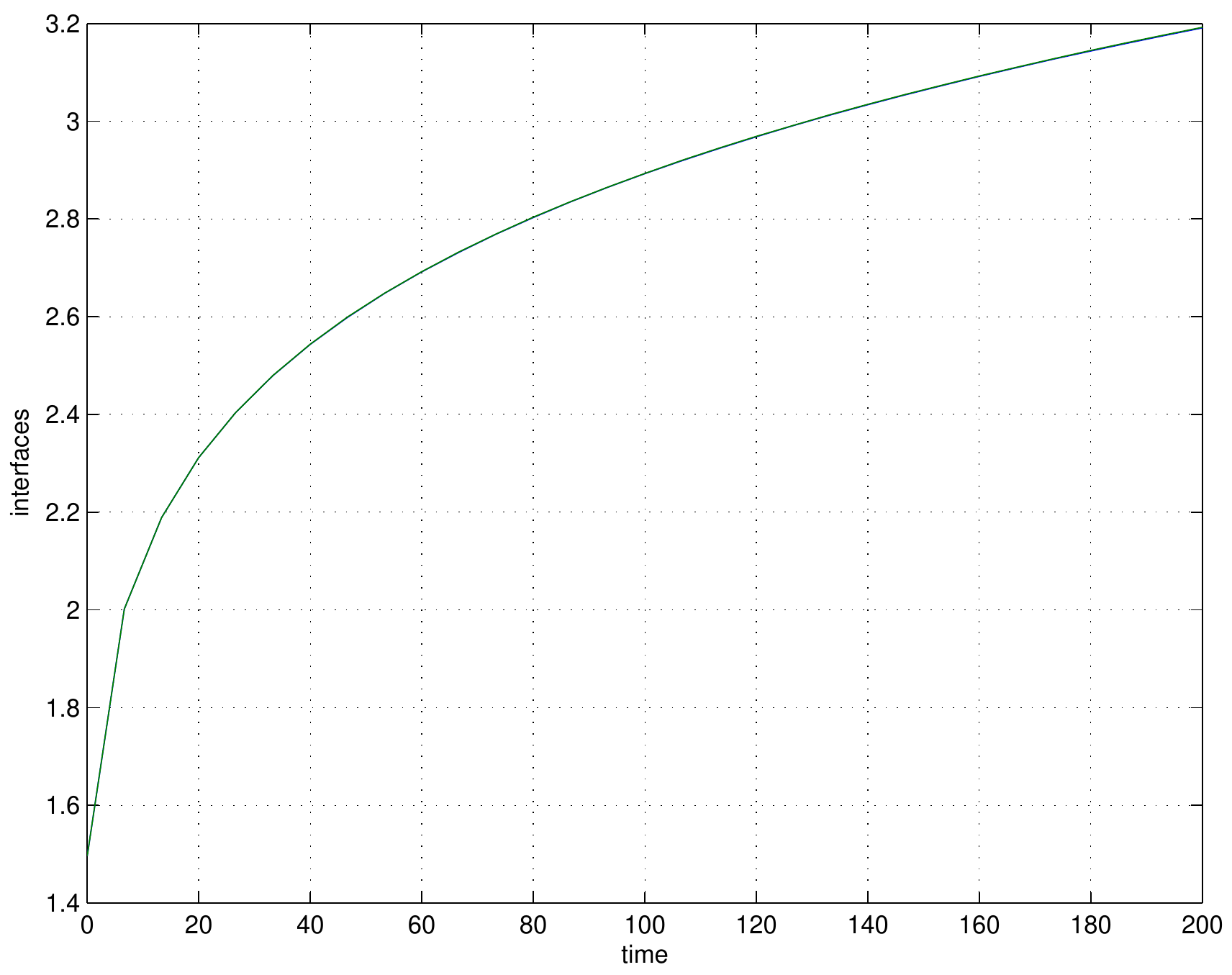}
    }}
\end{picture}
\caption{Time evolution of the interfaces (numerical and analytical)}    
\label{fig2}
\end{center}
\end{figure}

To compare the numerical $u_n$   and analytical $u_{an}$  solutions, we use a (relative) error estimator given by a numerical $L_2$-norm . We use discrete time points $t_j$ of the solutions approximated  in $x_i=s(t_j)y_i,\ i=1,...,N$ grid points and define
$$ 
L_{2,rel}(t_j)= \frac{\{ \sum_{i=1}^N |u_n(x_i,t_j)-u_{an}(x_i,t_j)|^2(y_i-y_{i-1})s(t_j)\}^{1/2}}{\{ \sum_{i=1}^N |u_{an}(x_i,t_j)|^2(y_i-y_{i-1})s(t_j)\}^{1/2}}$$
where the space integration is performed over the maximum length of supports of both solutions. In the same way we define the $L_{1,rel}(t_j)$ error (the exponent $2$ is replaced by $1$ in the integrands).
Mostly (in our experiments below) the time evolution of both errors does not differs significantly.
The time evolution of the relative $L_2$- error is shown in Fig. 3.
\begin{figure}
\begin{center}
  \setlength{\unitlength}{1cm}
\begin{picture}(12,8)
\put(0,0){\mbox{
    \includegraphics[width=12cm,height=6cm]{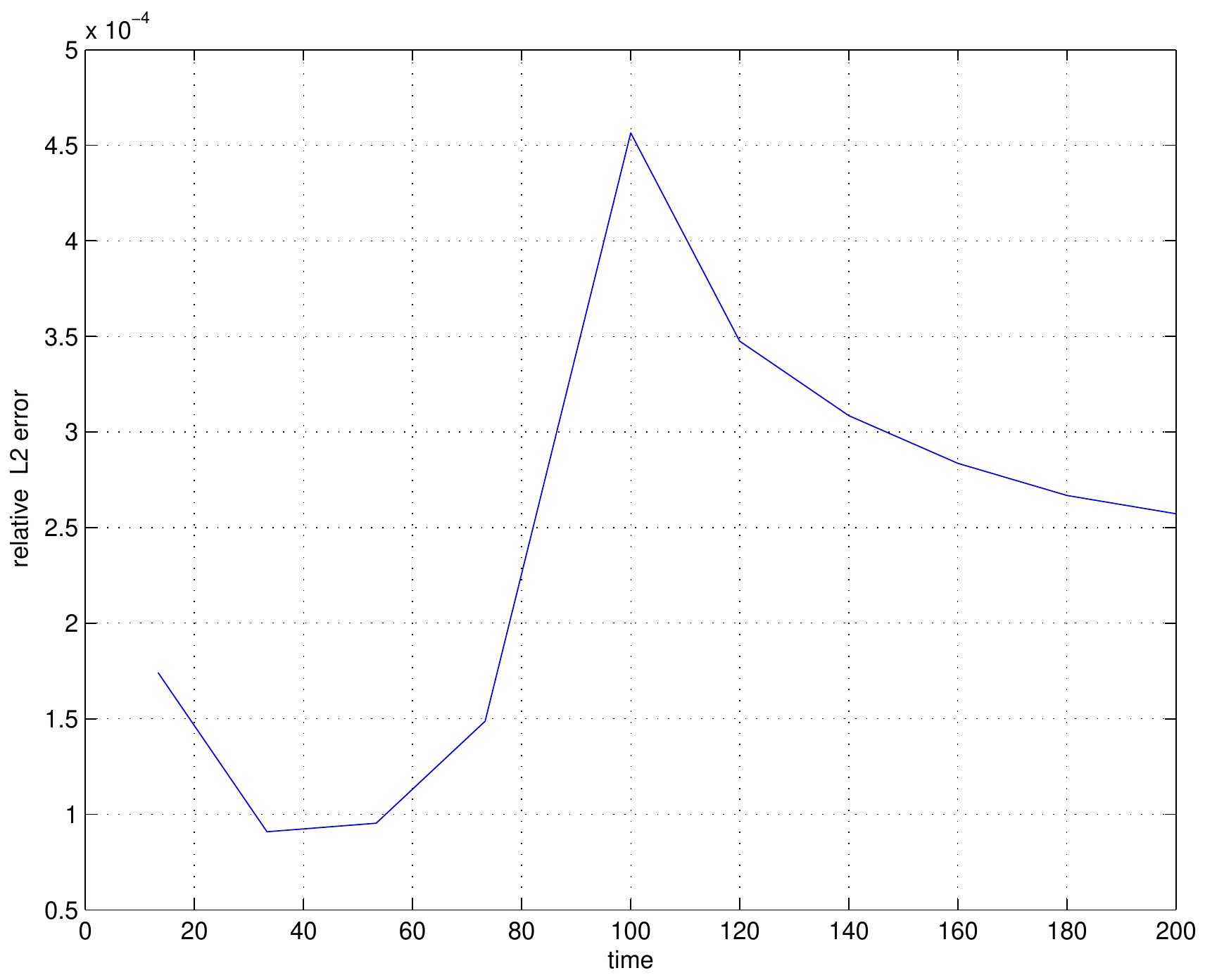}
    }}
\end{picture}
\caption{Time evolution of error}    
\label{fig3}
\end{center}
\end{figure}
Using more grid points ($N=40,M=15$),  the maximum $L_2$-error is $2.2 e-04$ and for ($N=60,M=20$) it is
 $1.3 e-04$.

Let us approximate (\ref{5}) in a similar way (reducing it to an ODE-system) without the interface modelling. We use a mass conservation scheme where
\begin{equation}
\label{kla}
\partial_x^2u(x_i,t)\approx \frac{2}{(\alpha(i)+\alpha(i+1) )(\alpha(i) \alpha(i+1))} 
   [\alpha(i+1)C_{i-1}.^p+\alpha(i)C_{i+1}^p-( \alpha(i+1)+\alpha(i))C_{i}^p)] =\dot{C}_i
\end{equation}
with $\alpha(i)=y_i-y_{i-1},\ i=1,...,N-1$. At the point $y=0$ we obtain the ODE
$$\dot{C}_0=\frac{2p}{\alpha(1)^2}(C_1-C_0)C_0^{p-1}$$
due to the symmetry arguments ($\partial_x u|_{x=0}=0$).
In Fig. 4 we present the corresponding numerical solution for $p=6$ with an equidistant space discretization of $x\in [0,10]$ with $N=M=100$ grid points. 
The error $ L_{2,rel}(t)\in (0.13,0.018)$ for $t\in (0,200)$. This error is nearly 100 times higher than the error of numerical results based on interface modelling with only $N=20$ grid points. The relative $ L_{2,rel}(t)$-error for the classical numerical approximation increases significantly with the reduction of grid points. The most of the contribution to the error is located near the interface.
\begin{figure}
\begin{center}
  \setlength{\unitlength}{1cm}
\begin{picture}(12,8)
\put(0,0){\mbox{
    \includegraphics[width=12cm,height=6cm]{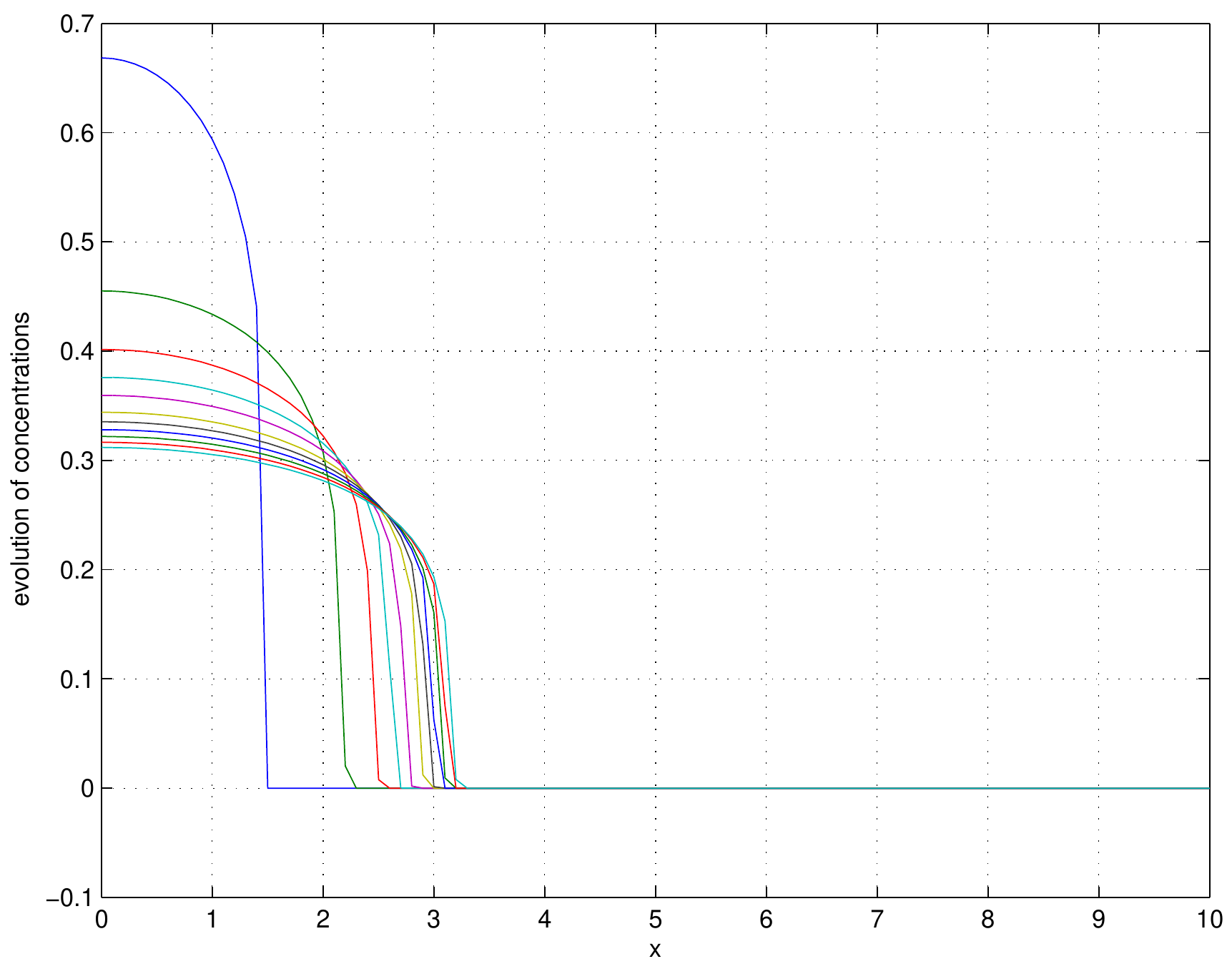}
     }}
 \end{picture}
\caption{Time evolution of the ``classical" numerical solution}    
\label{fig4}
\end{center}
\end{figure}

In the following Tables 1-3 we present the relation between the space discretization
$D_p,\ p=1,2,3,4$ and the average of $ L_{2,rel}(t)$-error versus the time 
$t\in (0,200)$ in terms of the following numerical approximation
$$ AL = \frac{1}{r} \sum_{j=1}^r  L_{2,rel}(t_j),$$
where we take $r=30$. We compare our numerical solution (based on interface 
modelling) with the analytical Barenblatt-Pattle solution for $n=6$.
In Table \ref{table1} we present $ N$ and $ AL$ where $M=N$ and discretization $D_4$
(i.e., unifom discretization in $y\in (0,1)$). In Table \ref{table2} we chose an optimal
$ 2\le M\le N$ so that $AL$ is minimal. In Table \ref{table3} we present
$N, D_p,d,M$ and $AL$ where for given $N$ and $D_p$ an optimal $M$ and $d$
are chosen so that $AL$ is minimal.
\begin{table}
\begin{minipage}[b]{0.48\linewidth}
\centering
\begin{tabular}{|l|r|}
\hline
\bfseries $N$ & \bfseries $AL. 10^4$\\ \hline 
10&2.66\\
20&2.45\\
30&2.35\\
40&2.1\\
50&2.2\\
100&2.3\\ \hline
\end{tabular}
\caption{$N$ versus $AL$ \label{table1}}
\end{minipage}
\hspace{0.5cm}
\begin{minipage}[b]{0.48\linewidth}
\centering
\begin{tabular}{|lc|r|}
\hline
\bfseries $N$ &  \bfseries M& \bfseries $AL. 10^4$\\ \hline 
10& 3& 1.85\\
20&5&1.2\\
30&5&0.67\\
40&7&0.585\\
50&10&0.57\\
100&30&0.51\\ \hline
\end{tabular}
\caption{Optimal $N,M$ versus $AL$ \label{table2}}
\end{minipage}
\end{table}
 
 \begin{table}
\begin{minipage}[b]{0.48\linewidth}
\centering
\begin{tabular}{|lccc|r|}
\hline
\bfseries $N$ &  \bfseries $D_p$&\bfseries $M$&\bfseries $d$& \bfseries $AL. 10^4$\\ \hline 
100&$D_3$&68& 3& 1.2\\
100&$D_2$&90& 3& 1.1\\
100&$D_1$&85& 6& 1.3\\
55&$D_3$&16&3&0.55\\
50&$D_2$&5&5&0.41\\
50&$D_1$&43&4&1.71\\ \hline
\end{tabular}
\caption{Optimal parameters versus $AL$\label{table3}}
\end{minipage}
\hspace{0.5cm}
\begin{minipage}[b]{0.48\linewidth}
\centering
\begin{tabular}{|l|r|}
\hline
\bfseries $N=M$ & \bfseries $AL. 10^3$\\ \hline 
200&4.3\\
100&9.6\\
50&22\\
30&36\\
20&55\\
10&112\\ \hline
\end{tabular}
\caption{$N=M$ versus $AL$ \label{table4}}
\end{minipage}
\end{table}

\subsection{An example of porous media equation with adsorption}

The time evolution of the interface can be very interesting when degenerated diffusion and adsorption is considered. In the following example we may see the dependence of the time evolution of the interface on the local profile of the solution (at the interface) and observe the order of the degeneracy very transparently. Let us consider the model problem

\begin{equation}
\label{Ke}
\partial_tu=\partial_x^2u^p-C_0u^{2-p}\quad \mbox{where} \ C_0>0,\    1<p<2
\end{equation}
with $\partial_x u_x|_{x=0}=0 $. For a special initial profile there exists an analytical solution
\cite{[RK]} given by the formula
$$ u(x,t)=a(t)^{-q_r}[\frac{C_0(p-1)^4\alpha^2+4p^2L_0^2}{4p^2((p-1)\alpha)^{2/(p+1)}}a(t)^{2/(p+1)}-$$
$$\frac{C_0(p-1)^2}{4p^2}a(t)^2-x^2]^{q_r},\quad \mbox{where}\ a(t)=\frac{2p(p+1)}{p-1}t+(p-1)\alpha , \  q_r=\frac{1}{p-1}$$
(when the value in parenthesis [.] is nonnegative, otherwise $u=0$) which corresponds to the initial profile
$$u_0(x)=\left\{
\begin{array}[a]{c c}
[(p-1)\alpha]^{-q_r}(L_0^2-x^2)^{q_r} &,\ \mbox{for}\ |x|<L_0\\
        0                      &, \ \mbox{for}\ |x|>L_0
\end{array}
\right. $$
where $\alpha>0,L_0>0$ are given parameters.\\
This example suits to our model considered in Theorem 1 (where $n=p,m=2-p, \forall p\in (1,2)$).
The interface is modelled (see Theorem 1) by
$$\dot{s}=-\frac{n}{n-1}\partial_xw+C_0(n-1)\frac{1}{\partial_xw}\quad \mbox{where}\ \beta =\frac{1}{p-1},\ u=w^{\beta},\ \beta=\frac{1}{p-1}.$$
Since we know an analytical solution (for the  special initial profile above) a very important question arises.
{\bf Does the analytical solution satisfy  our interface model given by Theorem 1?} The answer is positive and this can be verified by some additional computations. This supports  
 our results in interface modelling.
 In Fig. 5-8 we present the analytical solution, numerical solution, interfaces and relative error, respectively, for the parameters $p=1.8$, $N=100$ and $ M=80$. In this case the interface   starts to move forwards in the beginning  (there is a sufficiently high slope of the initial profile at the point $L_0=s(0)$) and  turns  backwards at the time $t\approx 5.5$, since then the influence of the adsorption prevails.

 \begin{figure}
\begin{center}
  \setlength{\unitlength}{1cm}
\begin{picture}(12,8)
\put(0,0){\mbox{
    \includegraphics[width=12cm,height=6cm]{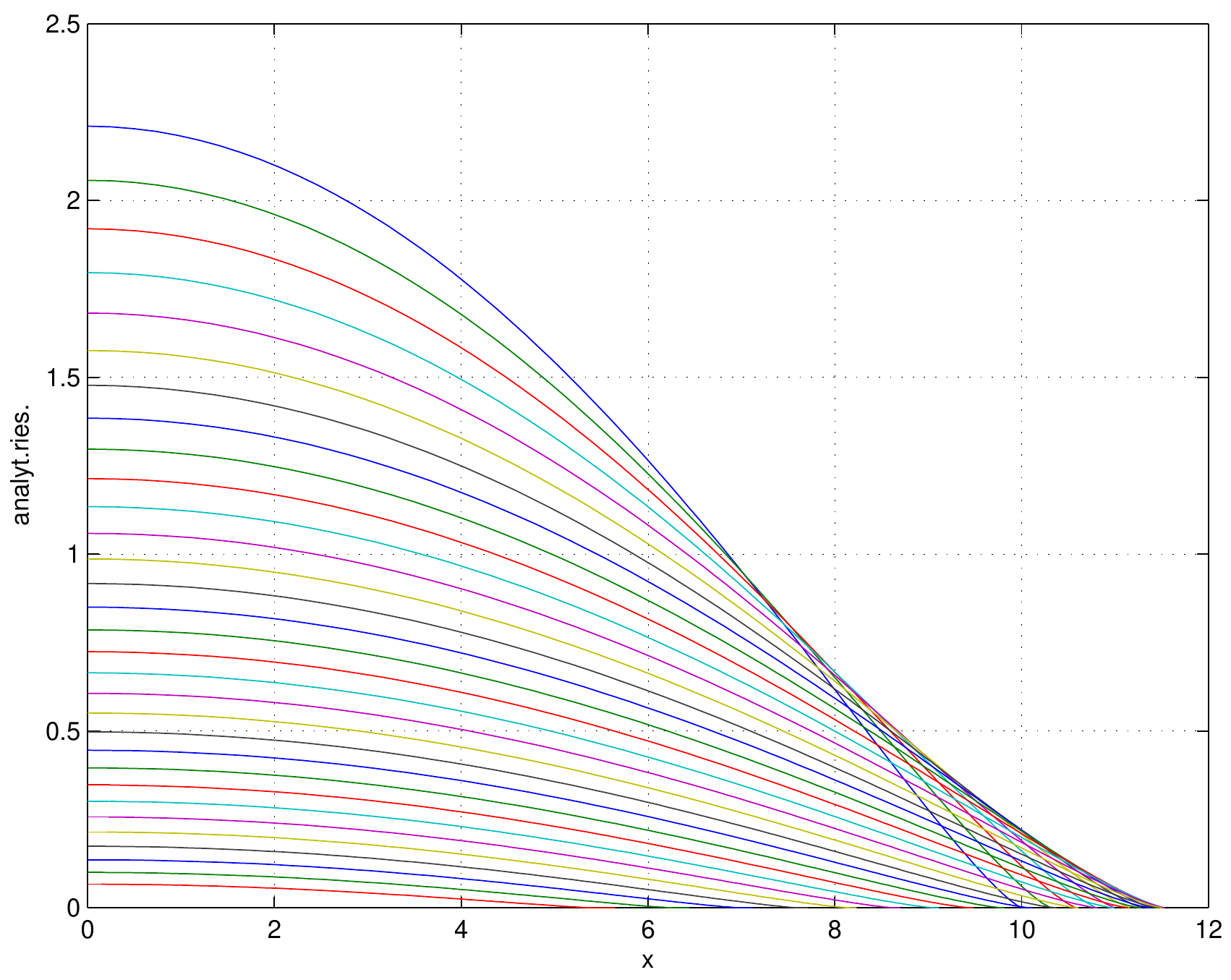}
     }}
 \end{picture}
\caption{The analytical solution in 31 equidistant time sections of $t\in (0,16)$}    
\label{fig5}
\end{center}
\end{figure}

 \begin{figure}
\begin{center}
  \setlength{\unitlength}{1cm}
\begin{picture}(12,8)
\put(0,0){\mbox{
    \includegraphics[width=12cm,height=6cm]{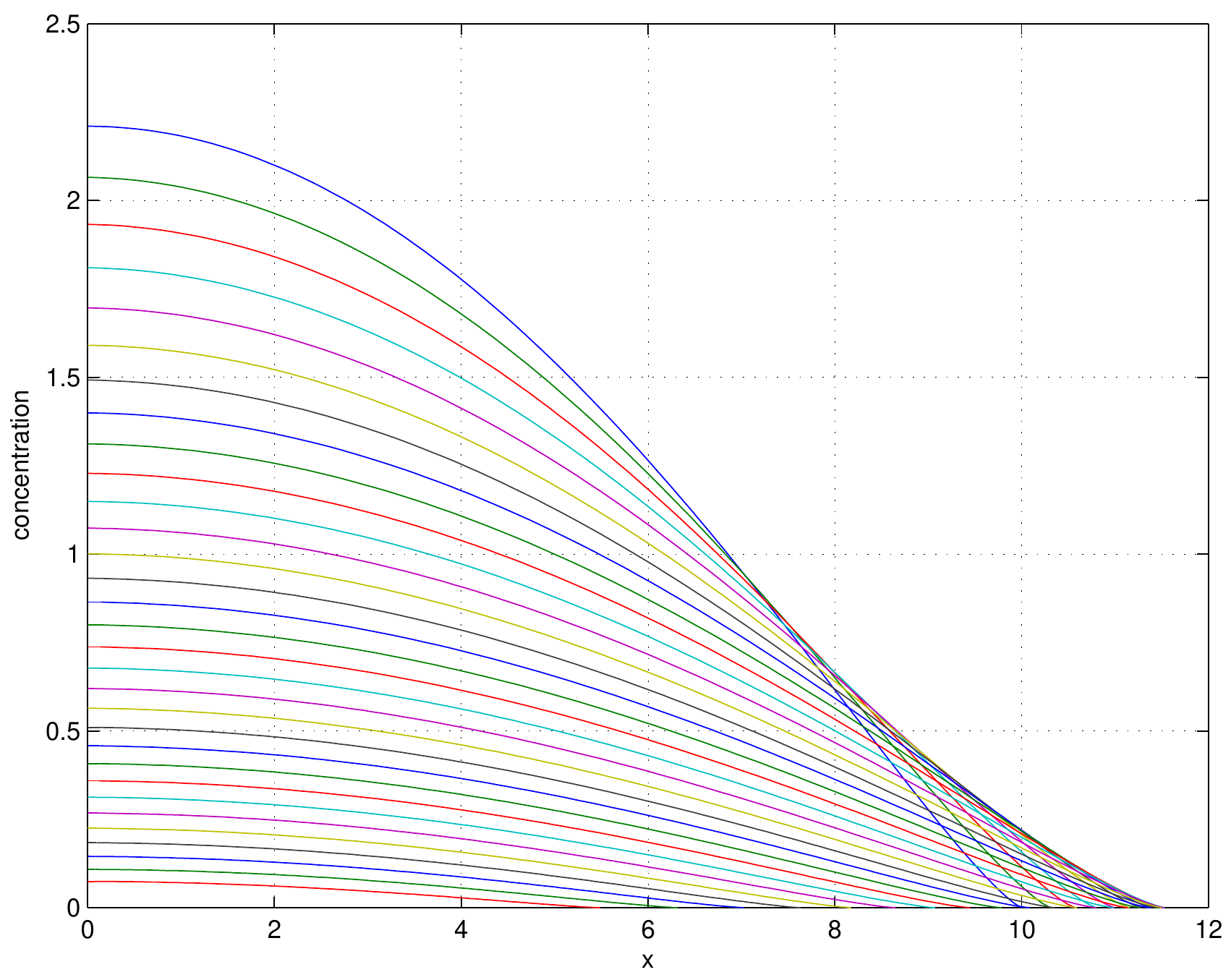}
     }}
 \end{picture}
\caption{The numerical solution in the same time sections as in Fig. 13, N=100}    
\label{fig6}
\end{center}
\end{figure}

 \begin{figure}
\begin{center}
  \setlength{\unitlength}{1cm}
\begin{picture}(12,8)
\put(0,0){\mbox{
    \includegraphics[width=12cm,height=6cm]{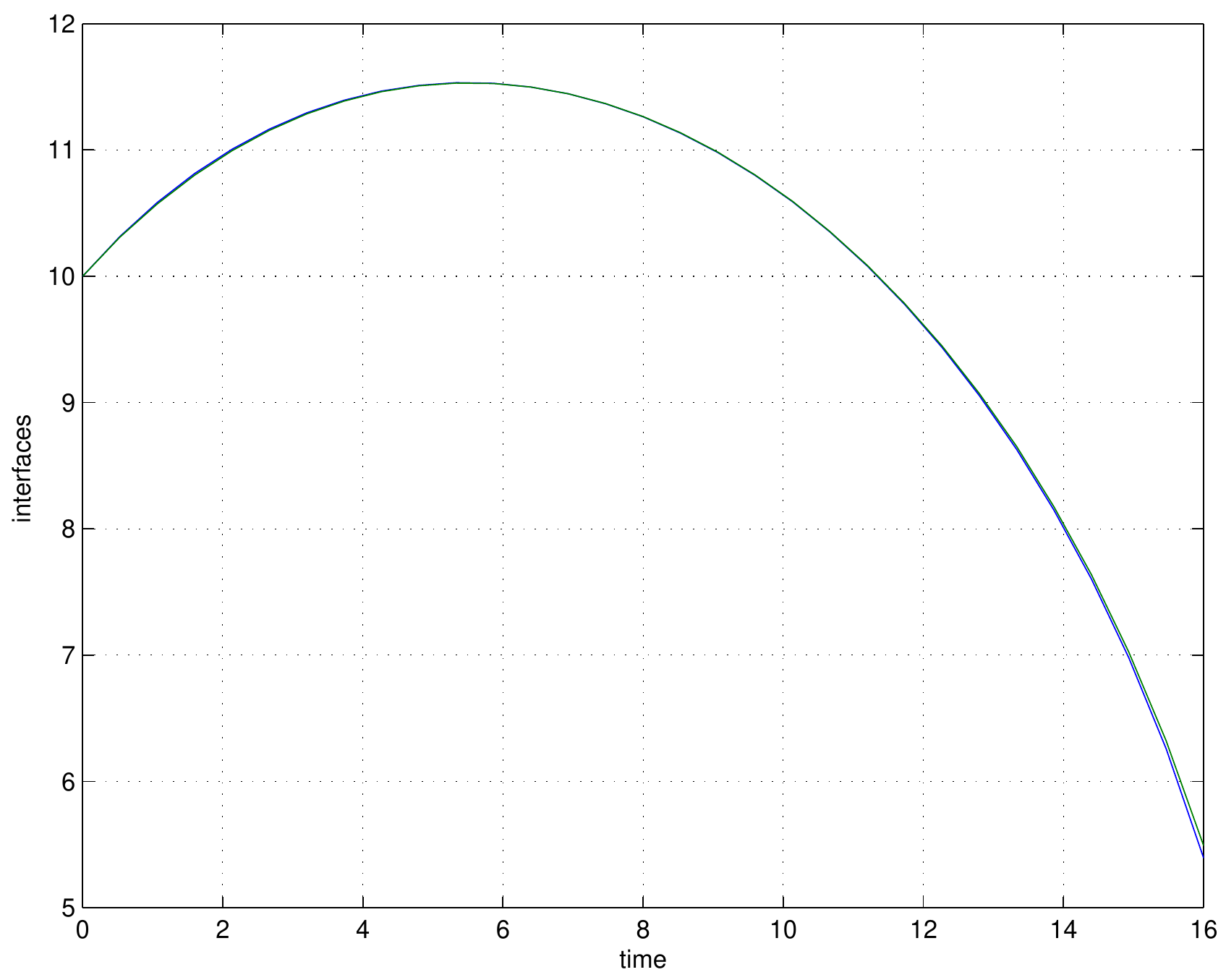}
     }}
 \end{picture}
\caption{Time evolution of interfaces (numerical and analytical), N=100}    
\label{fig7}
\end{center}
\end{figure}
 \begin{figure}
\begin{center}
  \setlength{\unitlength}{1cm}
\begin{picture}(12,8)
\put(0,0){\mbox{
    \includegraphics[width=12cm,height=6cm]{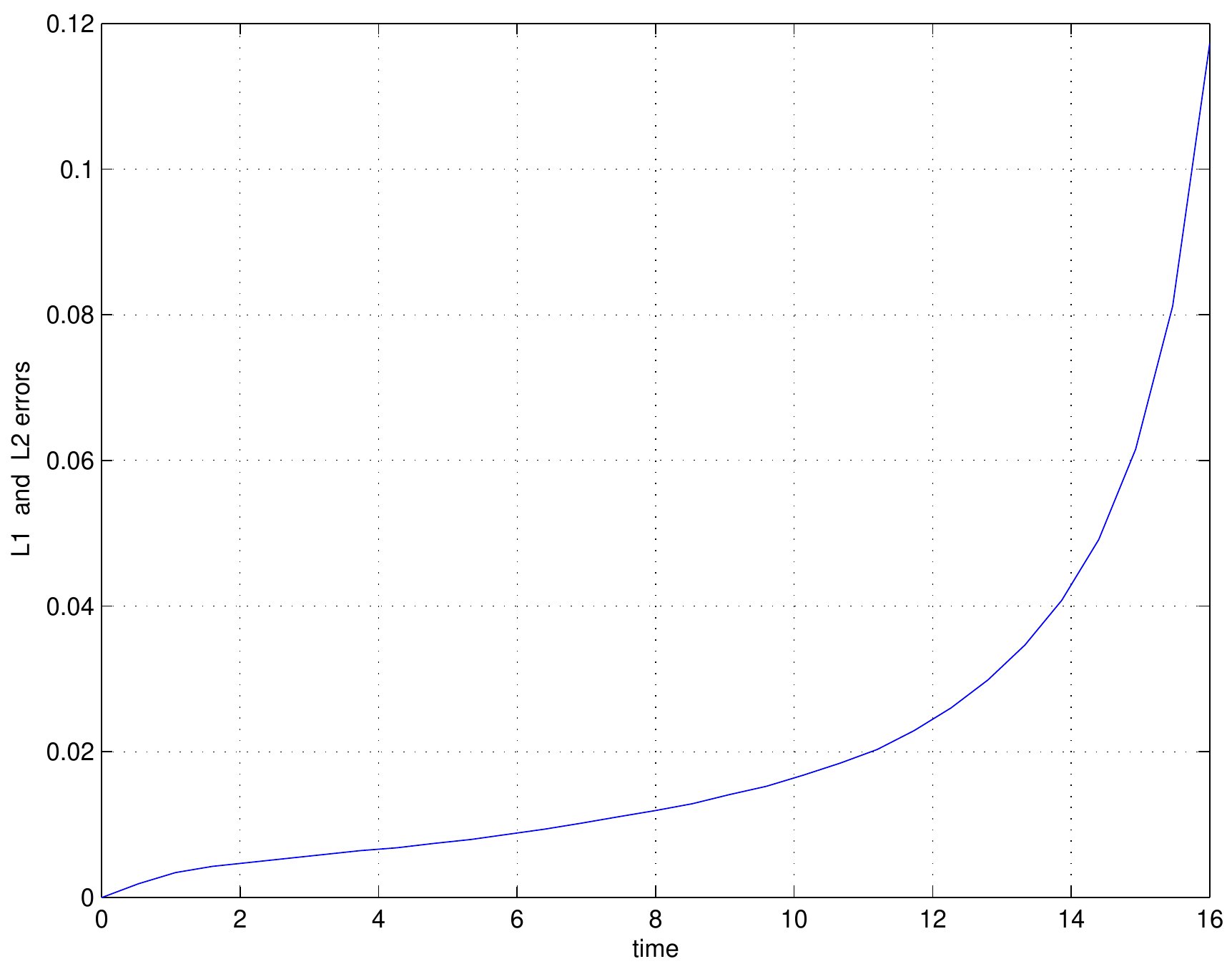}
     }}
 \end{picture}
\caption{Time evolution of errors, N=100 }    
\label{fig8}
\end{center}
\end{figure}
 At time $t_{que} \in (17,17.5)$ the solution is zero and this is a singularity point for our method. Also in the neighbourhood of $t=17$ the error significantly increases (the solution is very small). In Fig. 9  we present the interfaces  using only  $N=20$ grid points. In this case  $L_{2,rel}(t)\in (0,0.05)$ for $t\in (0,12.5)$ and then increases from $0.05$ up to $ 0.15$ for $t\in (12.5,17)$.
 \begin{figure}
\begin{center}
  \setlength{\unitlength}{1cm}
\begin{picture}(12,8)
\put(0,0){\mbox{
    \includegraphics[width=12cm,height=6cm]{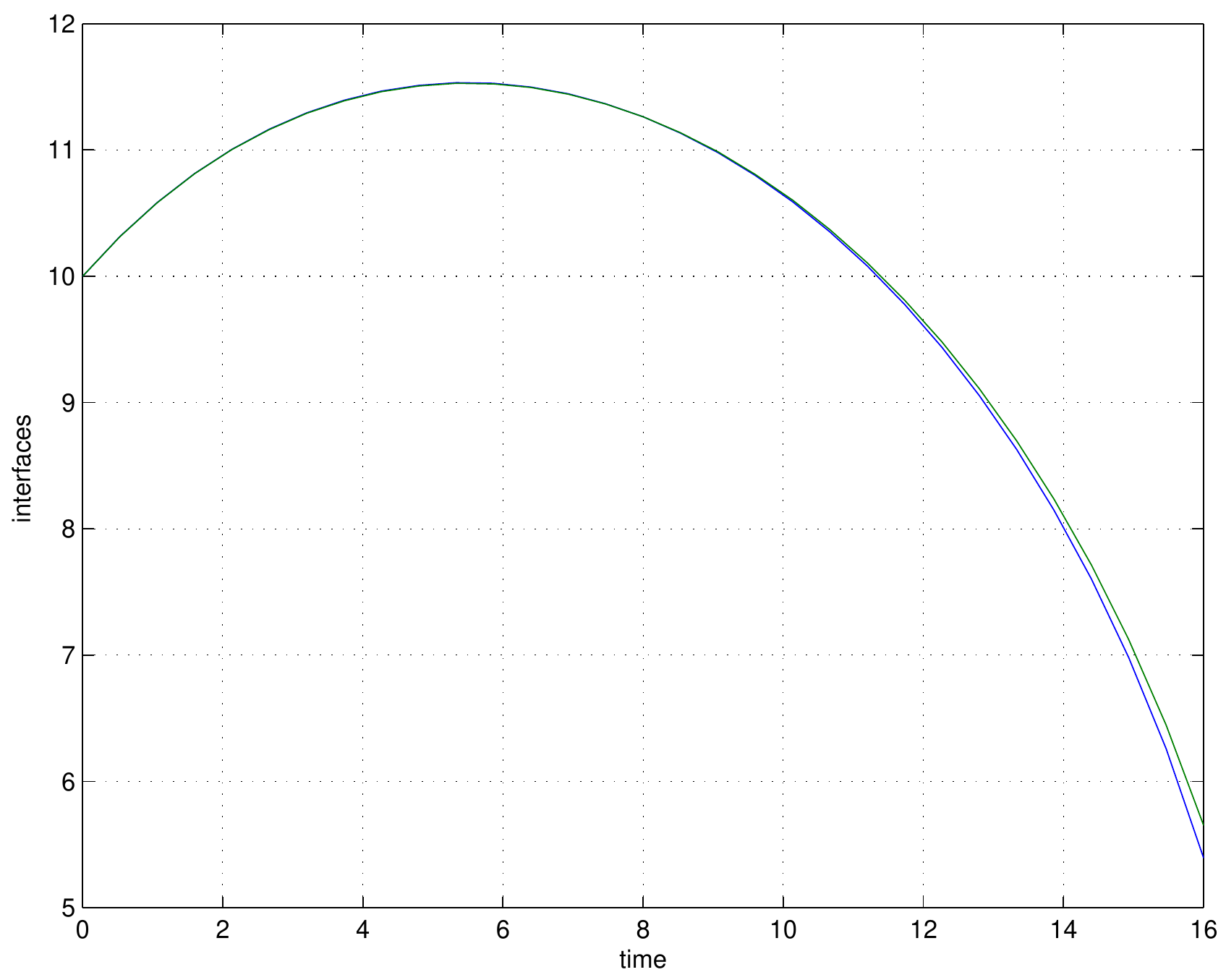}
     }}
 \end{picture}
\caption{Time evolution of interfaces (numerical and analytical),N=20}    
\label{fig9}
\end{center}
\end{figure}
As we can see even $N=20$ grid points are sufficient for the numerical approximation.\\
Similarly as in the previous example we solve this problem numerically without the interface modelling. In Fig. 10 we show the time evolution of the numerical solution  when using $N=100$ grid points. We consider the domain of the solution to be  $ (0,20)$.

 \begin{figure}
\begin{center}
  \setlength{\unitlength}{1cm}
\begin{picture}(12,8)
\put(0,0){\mbox{
    \includegraphics[width=12cm,height=6cm]{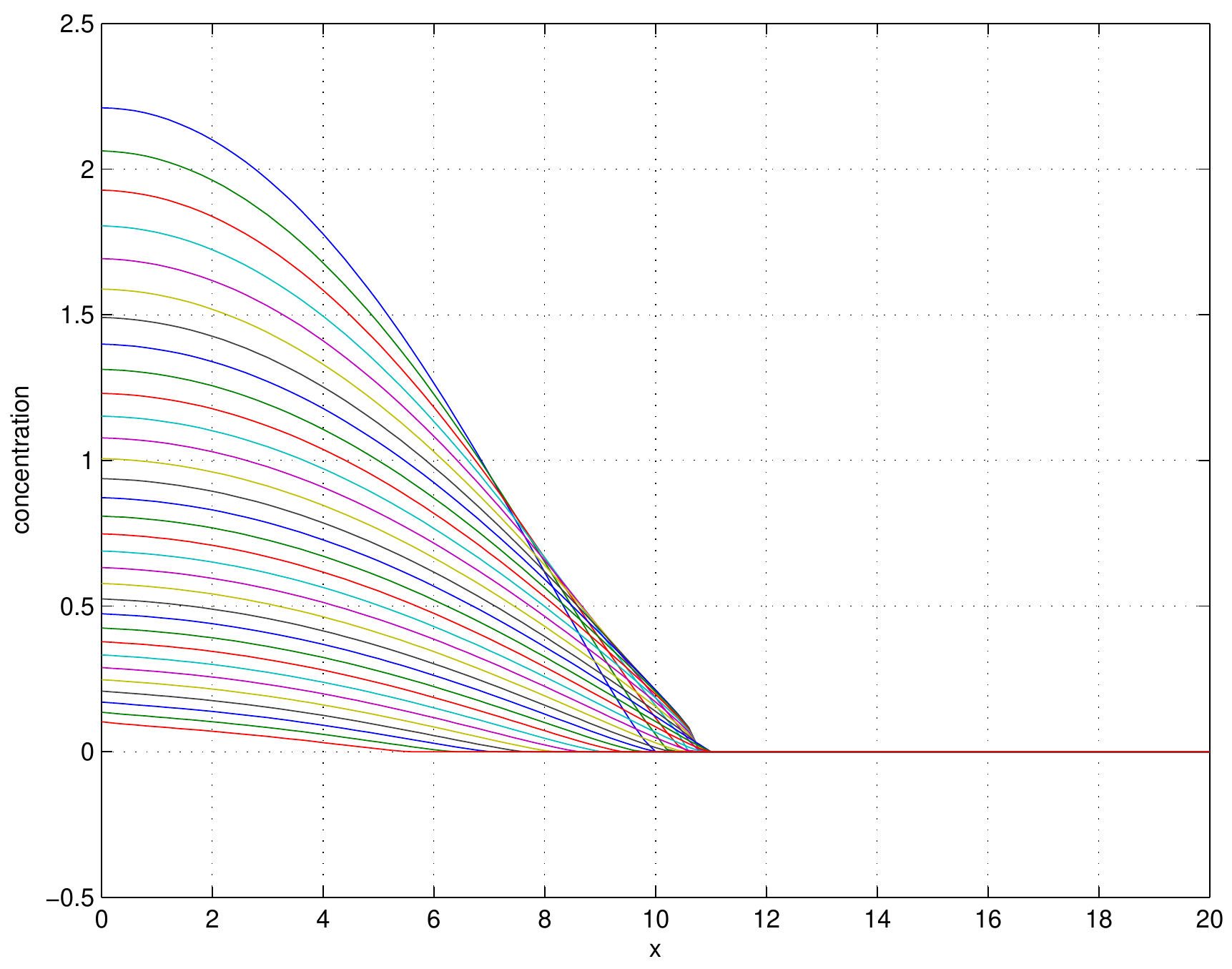}
     }}
 \end{picture}
\caption{Time evolution of the numerical solution, N=100}    
\label{fig10}
\end{center}
\end{figure}
In this case  $L_{2,rel}(t)\in (0,0.05)$ for $t\in (0,12)$ and then increasis from $0.05$ up to $ 0.25$ for $t\in (12,17)$.
When using only $N=40$ grid points, the time evolution of the error is 
  $L_{2,rel}(t)\in (0,0.5)$ for $t\in (0,8)$ and further increases from $0.5$ up to $5$ for $t\in (8,14)$.
The front of the solution is not as sharp as in the case of the Barenblat
solution. Here, the best discretization strategy is the uniform partition
($N=M $ of the interval $y\in (0,1)$) . The dependence of $AL$ on $N$ is presented in Table 4, where the time interval $(0,14)$ has been considered.\\

\subsection{Convection-diffusion-adsorption model\ (Contaminant transport)}

The governing equation of the (simplified 1D) contaminant transport  model with adsorption is of the form
\begin{equation}
\label{7}
\partial_tu =D\partial_x^2u -\partial_x(v u)-\rho \partial_tS,\ \mbox{in}\ x\in (0,L),t>0,
\end{equation}
coupled with the kinetics of adsorption
\begin{equation}
\label{8}
\partial_tS=\kappa(\Psi(u)-S),
\end{equation}
where $u$ is the concentration of the contaminant, $v$ is the velocity of the fluid in porous media
and $S$ is the adsorbed contaminant by the unit mass of the porous media with the density $\rho$. The function $\Psi$ is an adsorption isotherm, e.g., $\Psi(s)=c_0s^p, p\in (0,1)$ (Freundlich isotherm). When $\kappa \to \infty$ (adsorption is realized in an equilibrium mode), then $S=\Psi(u)$ and in this case (\ref{7}) is rewritten in the form
\begin{equation}
\label{9}
\partial_tB(u)=D\partial_x^2u-\partial_x(vu), \quad B(u)=u+\rho \Psi(u)
\end{equation}
which is of the form (\ref{4}) after the transformation $ \phi=B(u)$.\\
Numerical modelling of the interface (in the case of Freundlich isotherm) has been discussed in Section 5. For the construction of the interface model we do not invert $B$ and use (\ref{9})
$$\partial_tu = \frac{1}{B'(u)}(D\partial_x^2u-\partial_x(vu)$$
in the interface condition
$$ \dot{s}(t)=-\frac{\partial_tu}{\partial_xu}.$$
 In this case the interface is modelled by the ODE
$$ \dot{s}(t)=-\frac{D}{\rho a(1-p)} \partial_x w\quad  \mbox{for}\  x=s(t),\ \mbox{and}\ \ \beta=\frac{1}{1-p},\quad  u=w^{\beta}. $$
Implementing this interface model and rewritting (\ref{9}) in terms of $\phi$, we proceed as in Section 4 and obtain the following numerical results.
We use the parameters $p=1/2,\  t\in (0,1.9),\ N=150, M=40,\ \ v=1,\ D=0.05,\ a= \rho =1 $ and the Dirichlet condition $u(0,t)=1$  on the boundary $x=0$  . 
 The solution (concentration $u$) is presented in Fig. 11 in 11 equidistant time sections (starting from $t=0$), where the initial concentration is represented by the first (blue) graph (a regularization of the initial zero profile). The corresponding interface evolution is presented in Fig. 12. This will be denoted as the case I. Then, the final concentration profile $u(x,1.9)$ is taken as an initial profile for thecontaminant transport with the homogeneous
Dirichlet boundary condition $u(0,t)=0 ,\ \forall t \in (0,40)$ . The remaining parameters of the model are the same as in the previous experiment.  This will be denoted as the case II. The evolution of the concentration profile is presented in Fig. 13 and the corresponding interface in Fig. 14.

  \begin{figure}
\begin{center}
  \setlength{\unitlength}{1cm}
\begin{picture}(12,8)
\put(0,0){\mbox{
    \includegraphics[width=12cm,height=6cm]{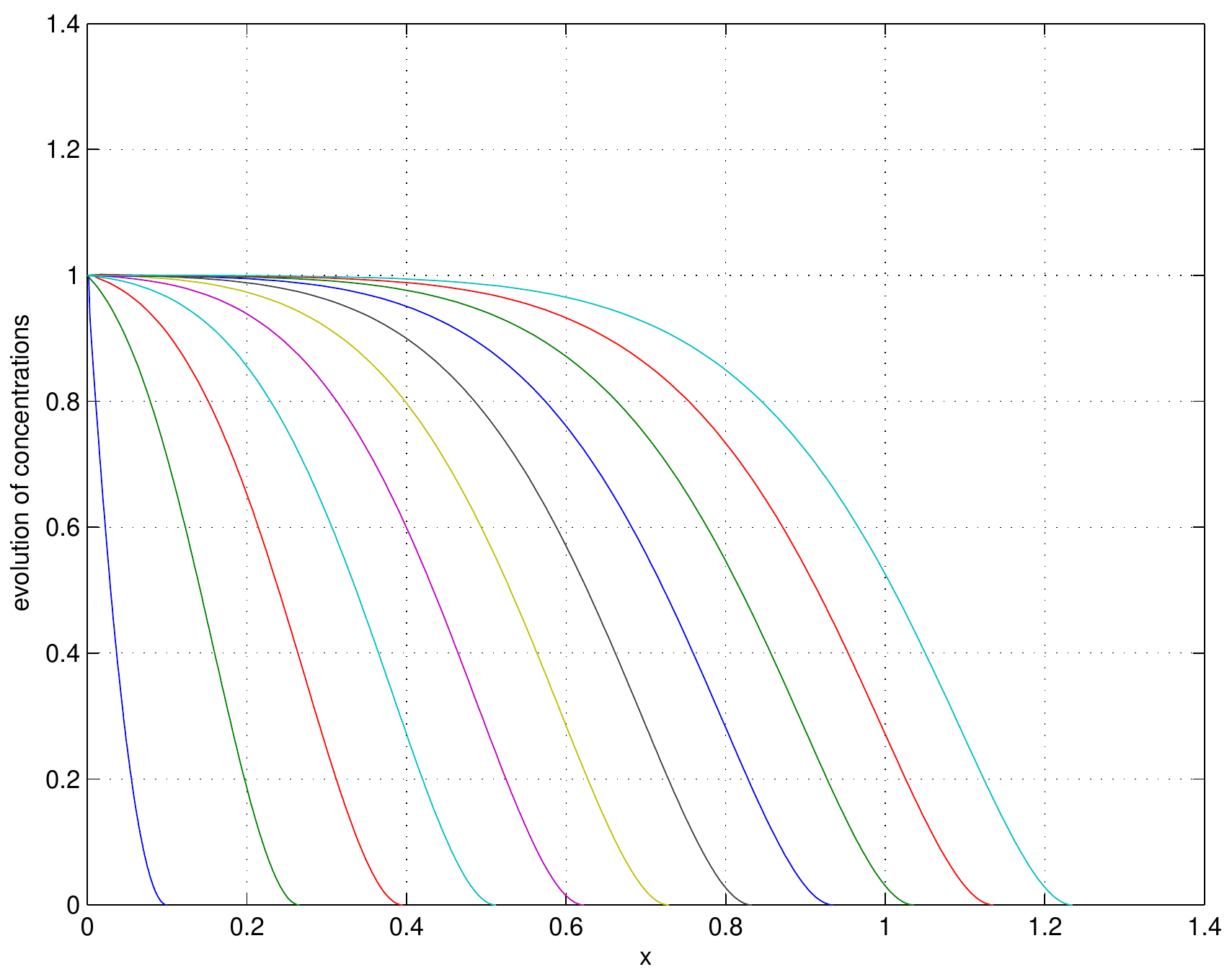}
     }}
\end{picture}
\caption{Time evolution of concentration I in 11 equidistant time sections,  p=0.5, N=150}    
\label{fig11}
\end{center}
\end{figure}

 \begin{figure}
\begin{center}
  \setlength{\unitlength}{1cm}
\begin{picture}(12,8)
\put(0,0){\mbox{
    \includegraphics[width=12cm,height=6cm]{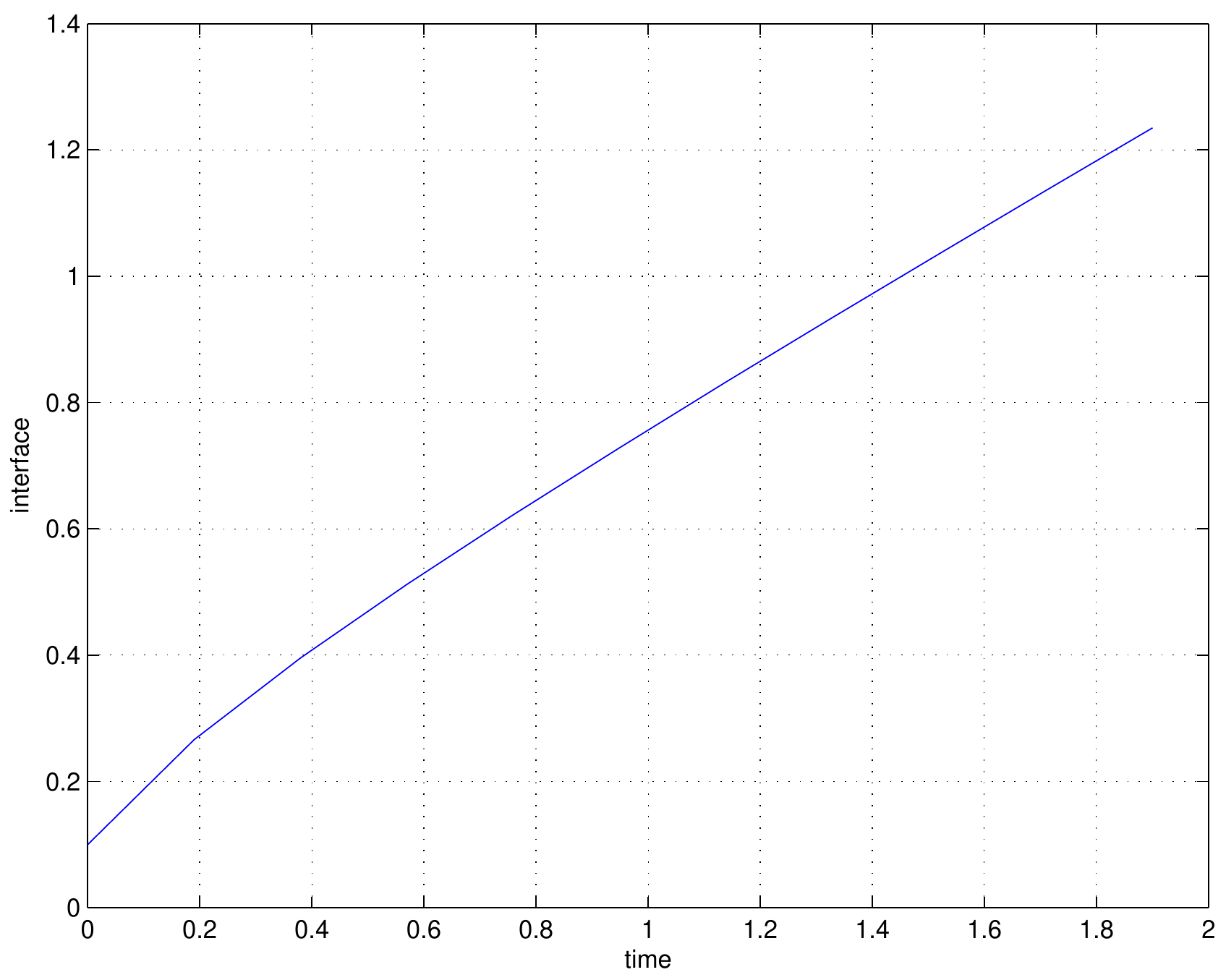}
     }}
 \end{picture}
\caption{Time evolution of interface I, p=0.5, N=150}    
\label{fig12}
\end{center}
\end{figure}

 \begin{figure}
\begin{center}
  \setlength{\unitlength}{1cm}
\begin{picture}(12,8)
\put(0,0){\mbox{
    \includegraphics[width=12cm,height=6cm]{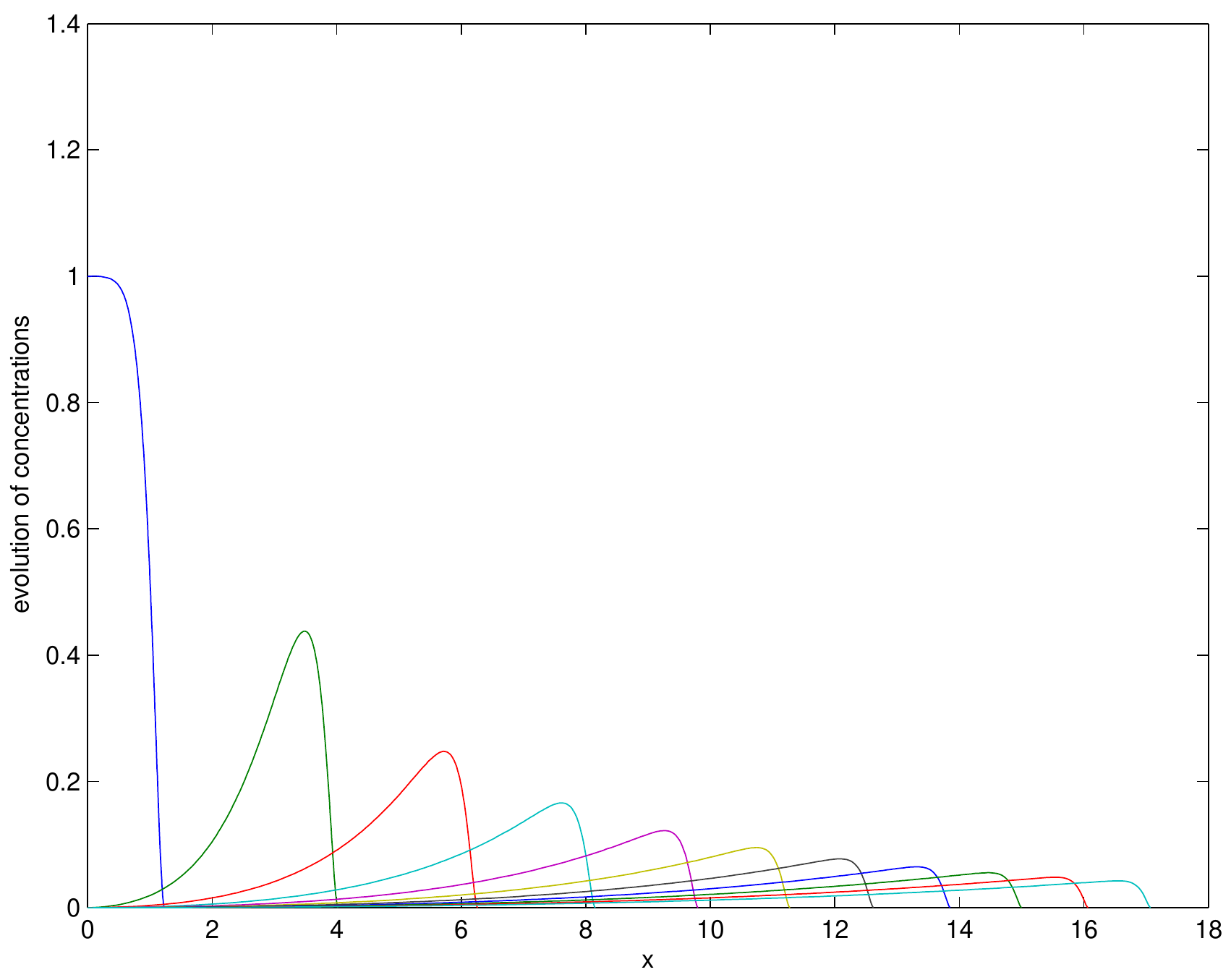}
     }}
 \end{picture}
\caption{Time evolution of concentrations II, p=0.5,\ M=60, N=150}    
\label{fig13}
\end{center}
\end{figure}

 \begin{figure}
\begin{center}
  \setlength{\unitlength}{1cm}
\begin{picture}(12,8)
\put(0,0){\mbox{
    \includegraphics[width=12cm,height=6cm]{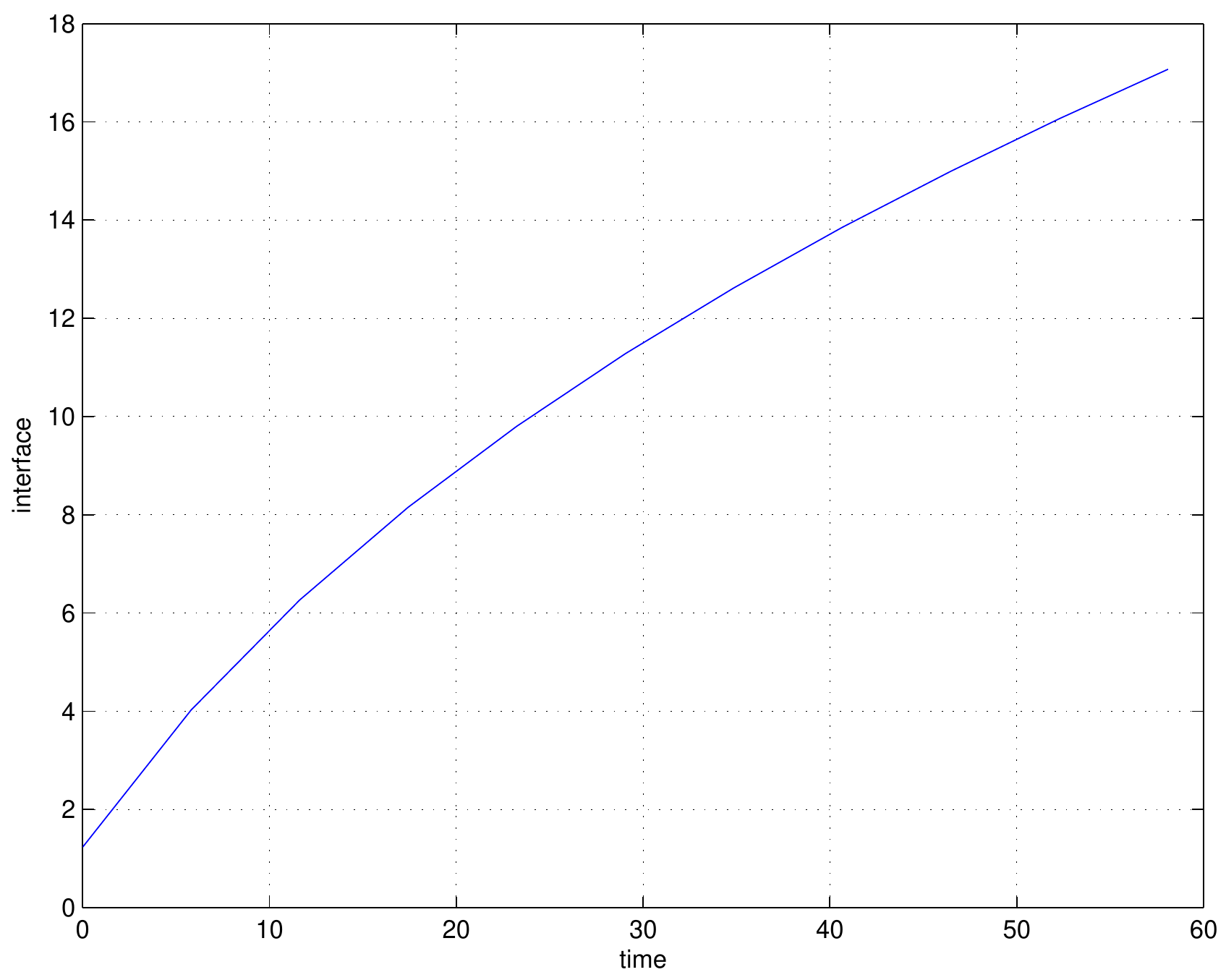}
     }}
 \end{picture}
\caption{Time evolution of interface II, p=0.5,M=60, N=150}    
\label{fig14}
\end{center}
\end{figure}
The numerical difficulty increases when we consider higher order degeneracy represented by $p=0.25$ and a small diffusion coefficient $D=0.005$. We keep other parameters from the previous experiment. In Fig. 15-16 we present the time evolution (  in $t\in (0,1.9)$-case I and $t\in (0,40)$-case II) of the concentration.

 \begin{figure}
\begin{center}
  \setlength{\unitlength}{1cm}
\begin{picture}(12,8)
\put(0,0){\mbox{
    \includegraphics[width=12cm,height=6cm] {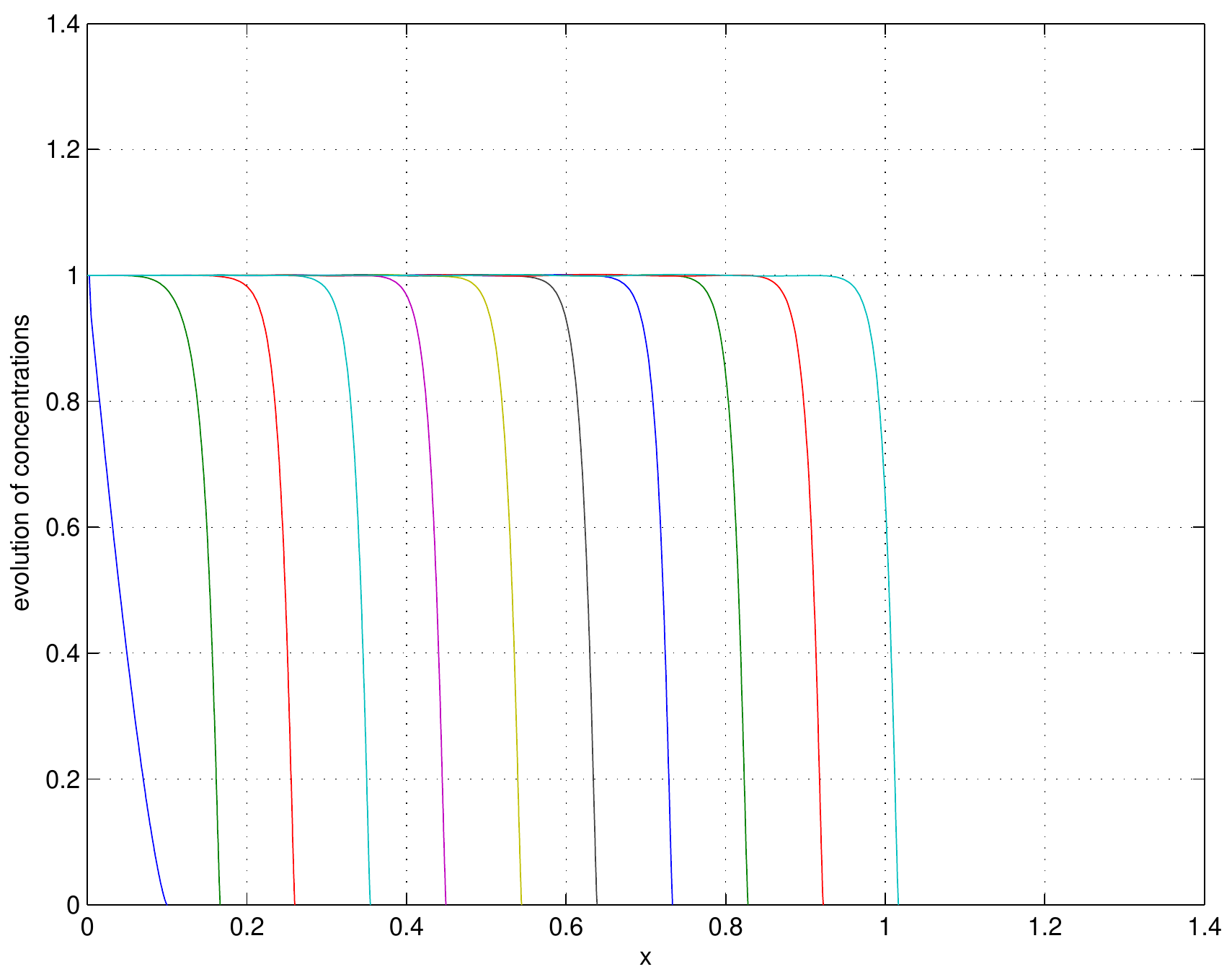}
     }}
 \end{picture}
\caption{Time evolution of concentration I, p=0.25, D=0.005, N=150}    
\label{fig15}
\end{center}
\end{figure}

 \begin{figure}
\begin{center}
  \setlength{\unitlength}{1cm}
\begin{picture}(12,8)
\put(0,0){\mbox{
    \includegraphics[width=12cm,height=6cm] {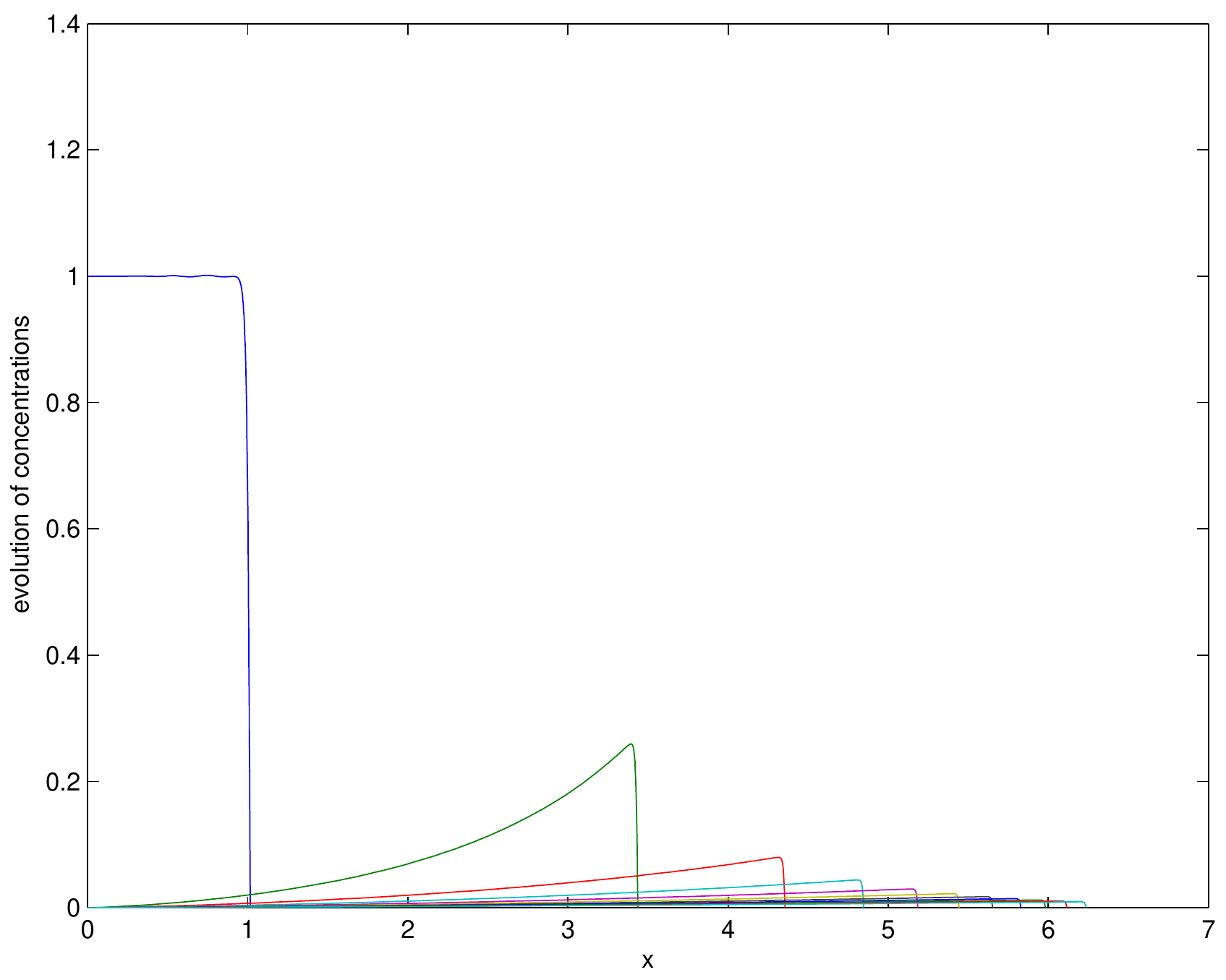}
     }}
 \end{picture}
\caption{Time evolution of concentration II, p=0.25, N=150, D=0.005}    
\label{fig16}
\end{center}
\end{figure}
The concentration profile drastically changes in Fig. 16. Therefore, we present the corresponding
solution in $11$ time sections of the interval $t\in (0,10)$ (in case II) - see Fig.
17. In this experiment the considered model is very near to
  nonlinear transport and the solution is very near to the entropy solution of the nonlinear transport with $D=0$
- see \cite{[FK]}. To increase the density of the grid points at the front we take $M=20,N=150$. 

 \begin{figure}
\begin{center}
  \setlength{\unitlength}{1cm}
\begin{picture}(12,8)
\put(0,0){\mbox{
    \includegraphics[width=12cm,height=6cm] {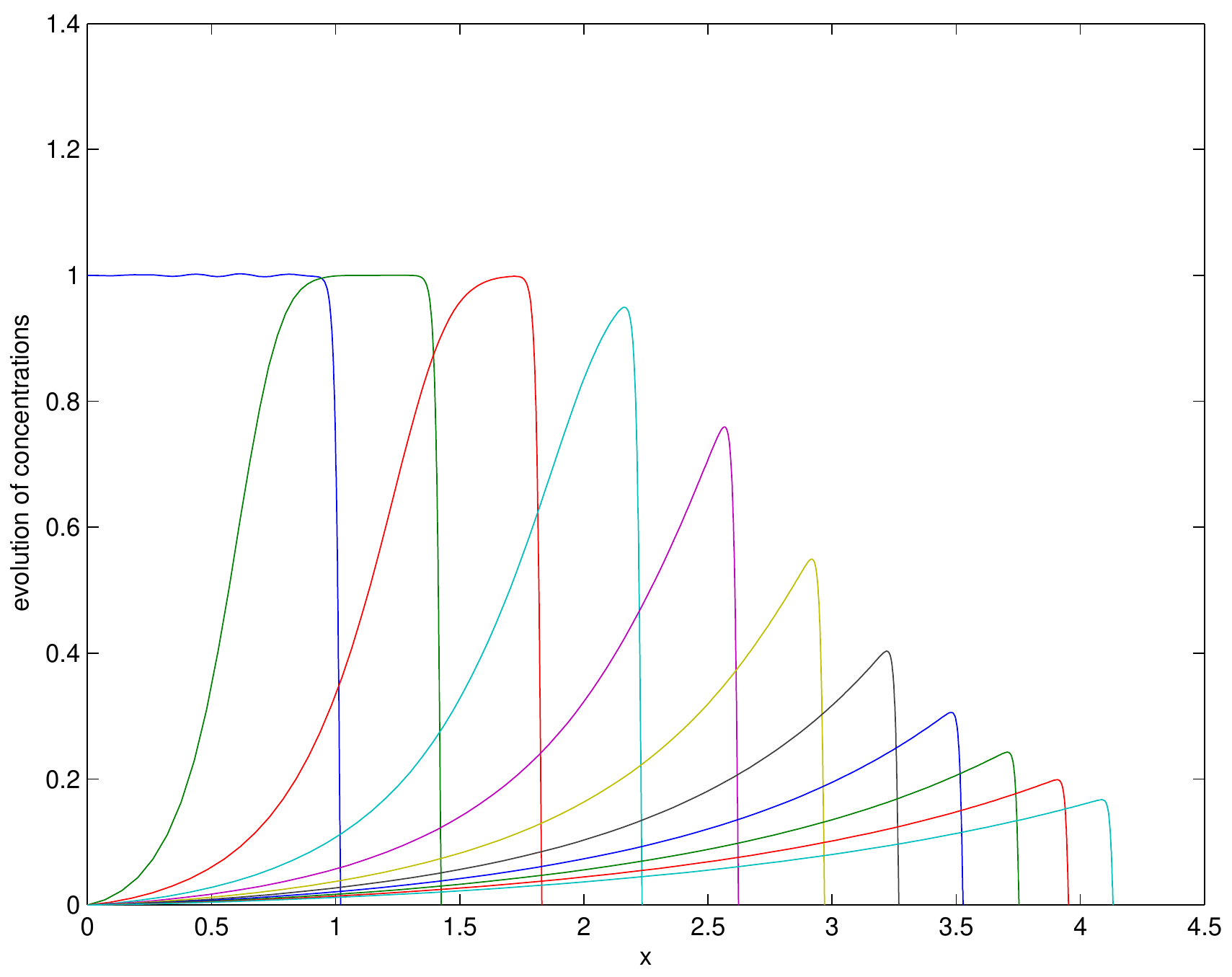}
     }}
 \end{picture}
\caption{Time evolution of concentrations II, $p=0.25,N=150, D=0.005, t\in (0,10)$}    
\label{fig17}
\end{center}
\end{figure}

In the fact the result presented in Fig. 17 demonstrates the efficiency of our numerical method.
The almost piecewise constant initial profile undergoes the nonlinear transport with zero boundary condition at $x=0$. This shock (at $x=0$) is not physically acceptable and immediately changes to the rarefaction which moves up to the front shock (physically acceptable) moving with the Rankin-Hugoniot speed. The solution endures also after the collision (rarefaction with the front shock). Now we are interested in finding out how many grid points are still sufficient for the numerical approximation.
In Fig. 18 we present the solution with the same parameters as in Fig. 17 with only $N=60$ ($ M=10$) grid points. The results are nearly the same as the time evolution of concentration with $N=200$ and $ M=20$ . This demonstrates the efficiency of our method.
 \begin{figure}
\begin{center}
  \setlength{\unitlength}{1cm}
\begin{picture}(12,8)
\put(0,0){\mbox{
    \includegraphics[width=12cm,height=6cm] {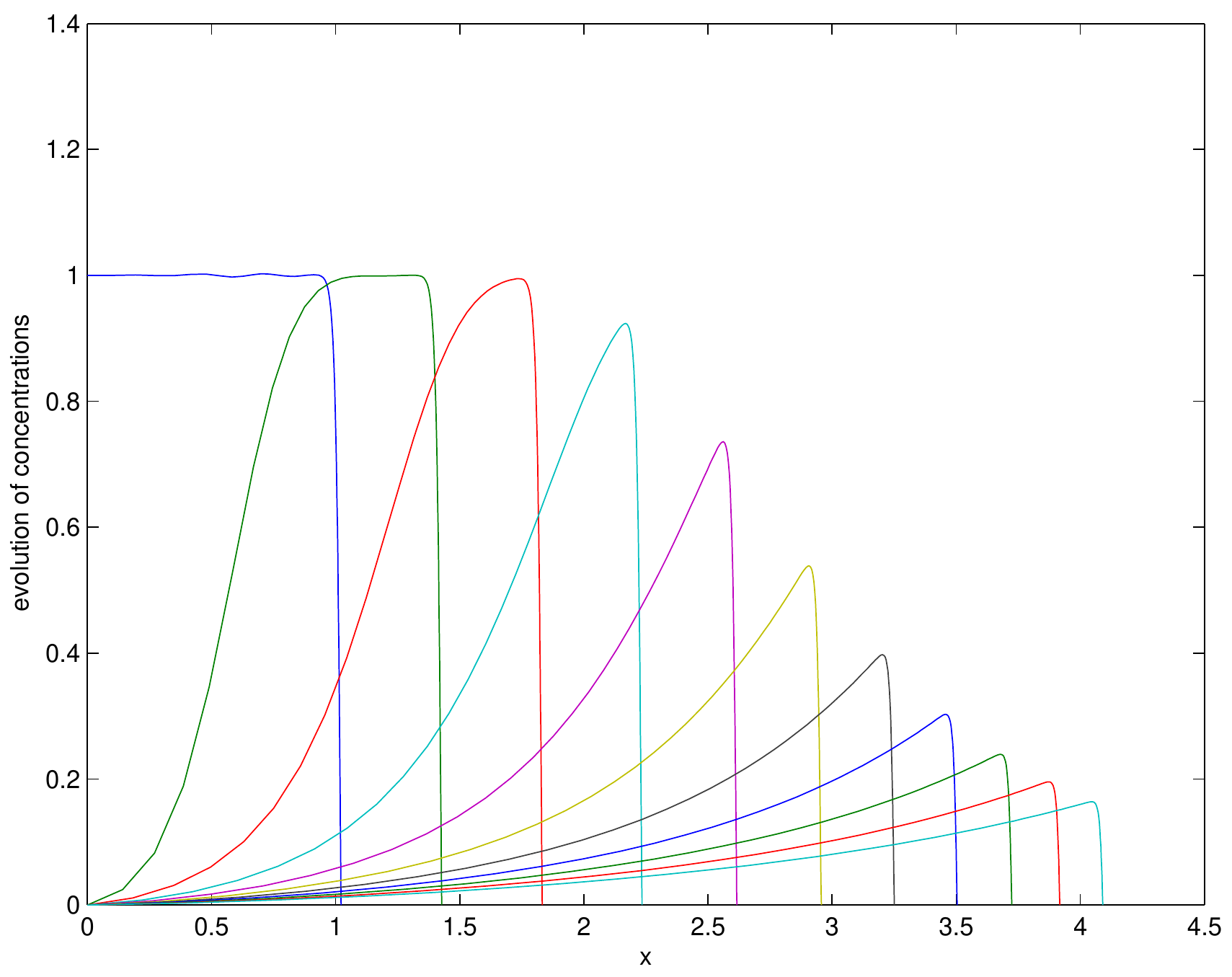}
     }}
 \end{picture}
\caption{Time evolution of concentrations $ p=0.25, N=60, D=0.005, t\in (0,10)$}    
\label{fig18}
\end{center}
\end{figure}
The experiments with higher degeneracy $p<0.85$ are more suitable
 for our numerical approximations. But when $p\to 1$ ($p=1$ is linear case) troubles arise with our numerical approximations since this is  the singular case for our method.\\
The same mathematical model in 1D was  treated numerically in \cite{[FK]} using the operator splitting method, where in the transport part a semi-analytical solution has been used and in the diffusion part the finite volume method has been applied. A very precise numerical solution of this model have been applied in a 2D-problem for the transport of contaminant in a dual-well setting. The jump in the Dirichlet boundary condition is applied there, to create the pulse shape in the injection well and to compute (and also to measure) the response in the extraction well. This gives very important information for the solution of the inverse problems (in the determination of the model parameters). Since it is an ill-posed problem, a very precise numerical solution is required. Therefore our method can be successfully applied there and we will implement it in our forthcoming paper. The method used in \cite{[FK]} is limited to Freundlich and Langmuir isotherms. Our method also works in a much more general case of sorption isotherms, since it does not depend on their geometrical properties. For example,
 we consider $\Psi(s)=\frac{a s^p}{1+bs^p}$ with $p\in (0,1)$. This isotherm equals to the Freundlich one
 for $b=0$ and to the Langmuir one  for $p=1$ . In this setting our model (\ref{9}) admits the solution with interface ($p\in (0,1)$). The case $p=1$ is reduced to the Langmuir sorption isotherm which is a nondegenerated problem. This case is a singular one for our method (the speed of interface is $\infty$). However, we take $p=0.95$ and expect that our solution will be close to the one with Langmuir sorption isotherm solved in \cite{[FK]} with $D=0$ (only transport). To compare these results, we take $ D=0.005$, $a=1.5$ and $b=1$. We use $N=100$ and $ M=20$ in $D_4$ discretization. The corresponding results are in Fig. 19, which we can compare (graphically) with the corresponding results in \cite{[FK]} (Fig. 4). Our concentration profiles at time sections 5 and 11 then correspond to the ones for $t=16$ and $ t=40$ in \cite{[FK]}. There is a good agreement. This gives us the possibility to 
 obtain a good approximation also in the case of sorption isotherms $\Psi(s),\  |\dot{\Psi)}'(0)|<\infty $ representing a nondegenerated problem (without interface). In that case, in the place of $\Psi(s)$ we consider an approximation $\Psi(s^p)$ with $p\nearrow 1$, which is of course limited by the numerical stability of the ODE solver in our setting. We observe that the $p=0.95$ leads to a good approximation.
 
  \begin{figure}
\begin{center}
  \setlength{\unitlength}{1cm}
\begin{picture}(12,8)
\put(0,0){\mbox{
    \includegraphics[width=12cm,height=6cm] {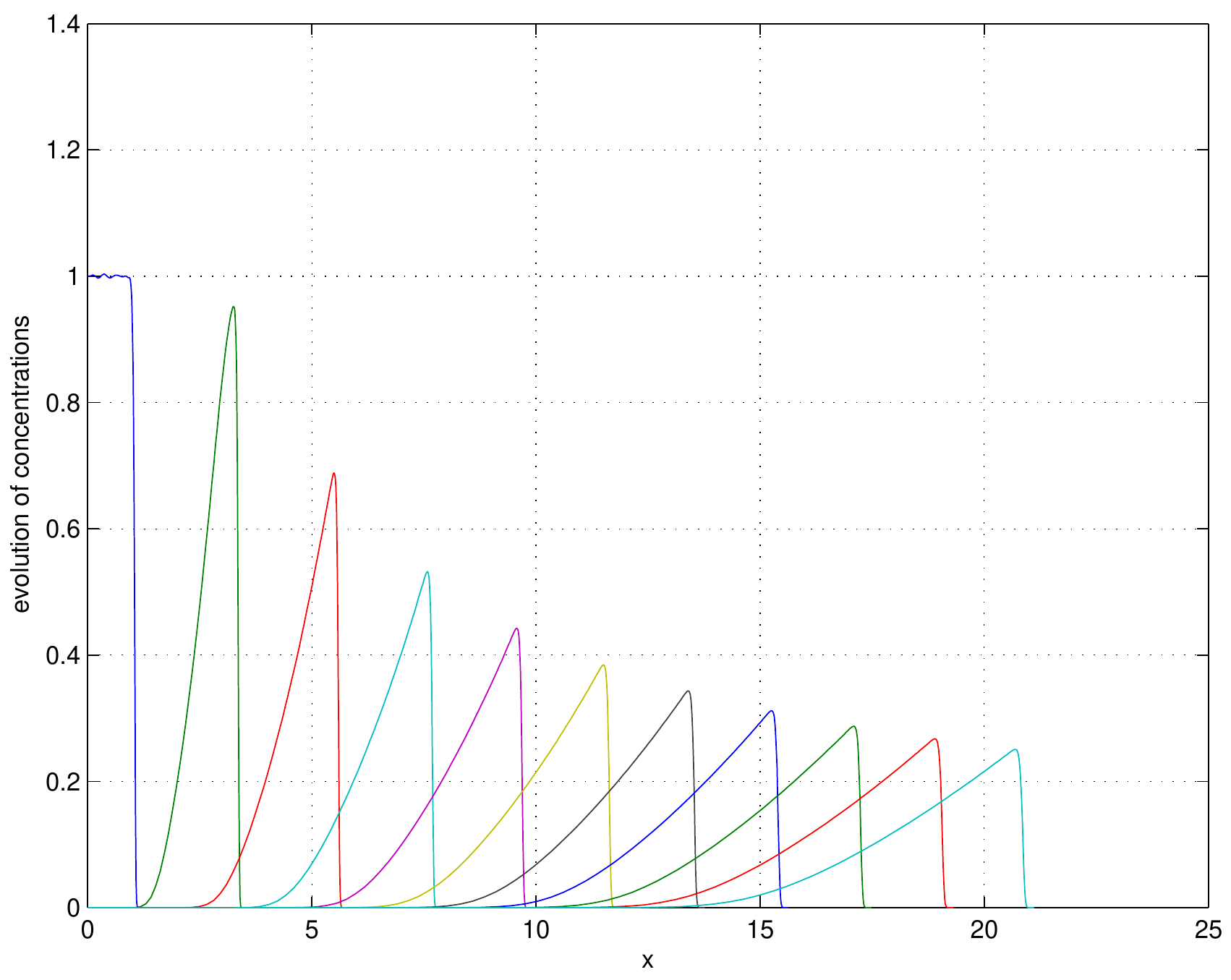}
     }}
 \end{picture}
\caption{Time evolution of concentrations $a=1.5, b=1 p=0.95, N=100, D=0.005, t\in (0,40)$}    
\label{fig19}
\end{center}
\end{figure}
In the next experiment we take $b=0.01$ and the remaining parameters are the same as  in the previous experiment (with $p=0.95$).
The results  are shown in Fig. 20 and they represent a good approximation of linear sorption (case $p=1$). Consequently, when a sorption isotherms $\Psi(s)$ doesn't leads to the appearance of the interface (nondegenerated case), we can use the approximation based on the Langmuir sorption type $\Psi(\frac{s^p}{1+bs^p})$ with $p\nearrow 1,\ b\searrow 0$ close to the instability region.
  \begin{figure}
\begin{center}
  \setlength{\unitlength}{1cm}
\begin{picture}(12,8)
\put(0,0){\mbox{
    \includegraphics[width=12cm,height=6cm] {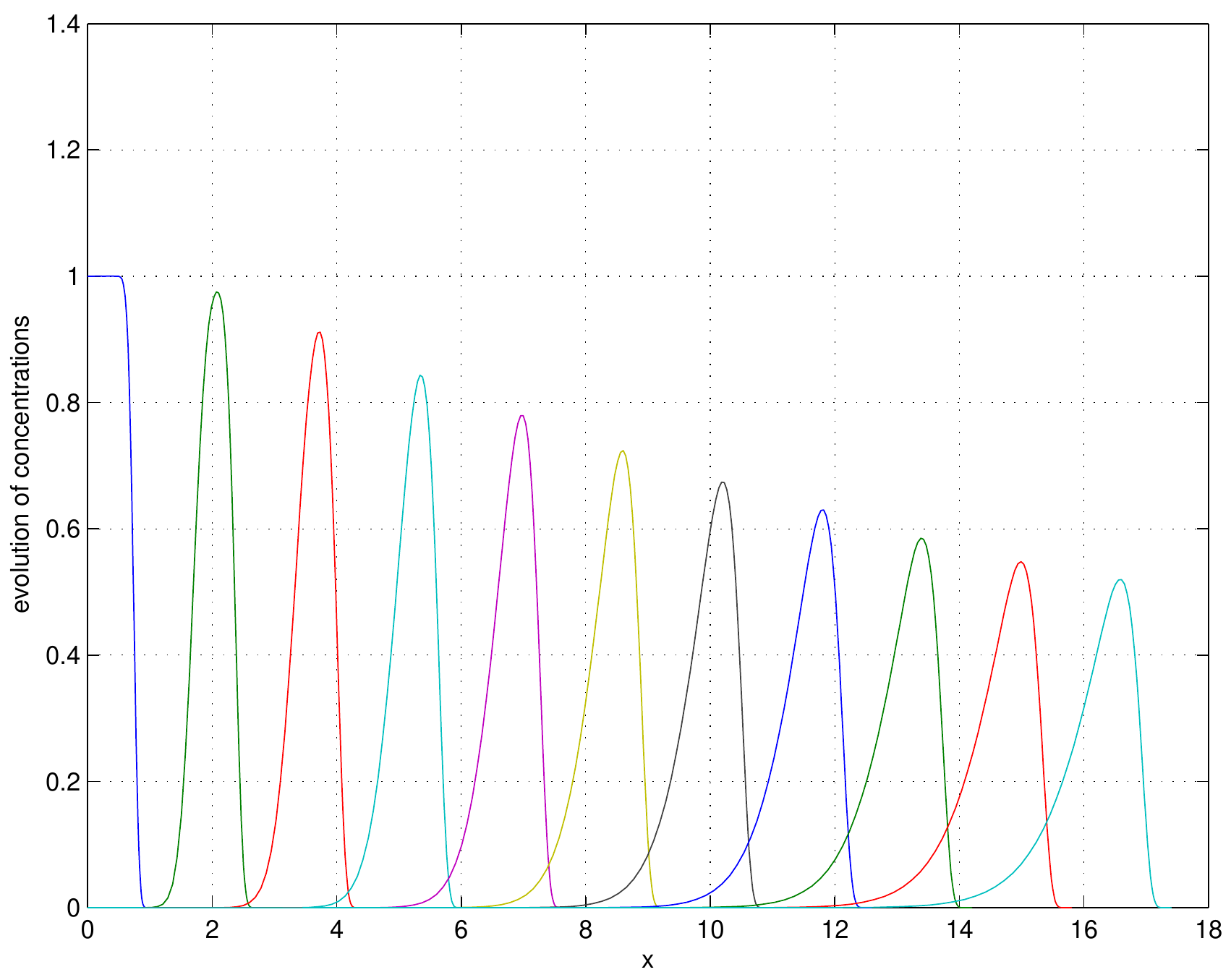}
     }}
 \end{picture}
\caption{Time evolution of concentrations $a=1.5, b=0.01, p=0.95, N=100, D=0.005, t\in (0,40)$}    
\label{fig20}
\end{center}
\end{figure}

Now we compare our results  with the corresponding ones in \cite{[R1]} (based on the operator splitting method) with the following parameters: $ p=0.75, D=0.1, V=1$ and $ a=2$ . We will use the discretization parameters $ N=100$ and $ M=20 $.
Time evolutions of the concentrations are presented in Fig 21. The graph corresponding to the $5$-th time section
 can compared with the one in \cite{[R1]} (Fig. 1 for $t=15$). Then, we  use $D=0.001$ and other parameters are the same as in the previous experiment. Then, we can compare the graph of the solution shown in Fig. 22 ($5$-th time section) with the corresponding results in  \cite{[R1]} ( Fig. 2  for $t=15$ ). In both cases we obtain a very good agreement.
  \begin{figure}
\begin{center}
  \setlength{\unitlength}{1cm}
\begin{picture}(12,8)
\put(0,0){\mbox{
    \includegraphics[width=12cm,height=6cm] {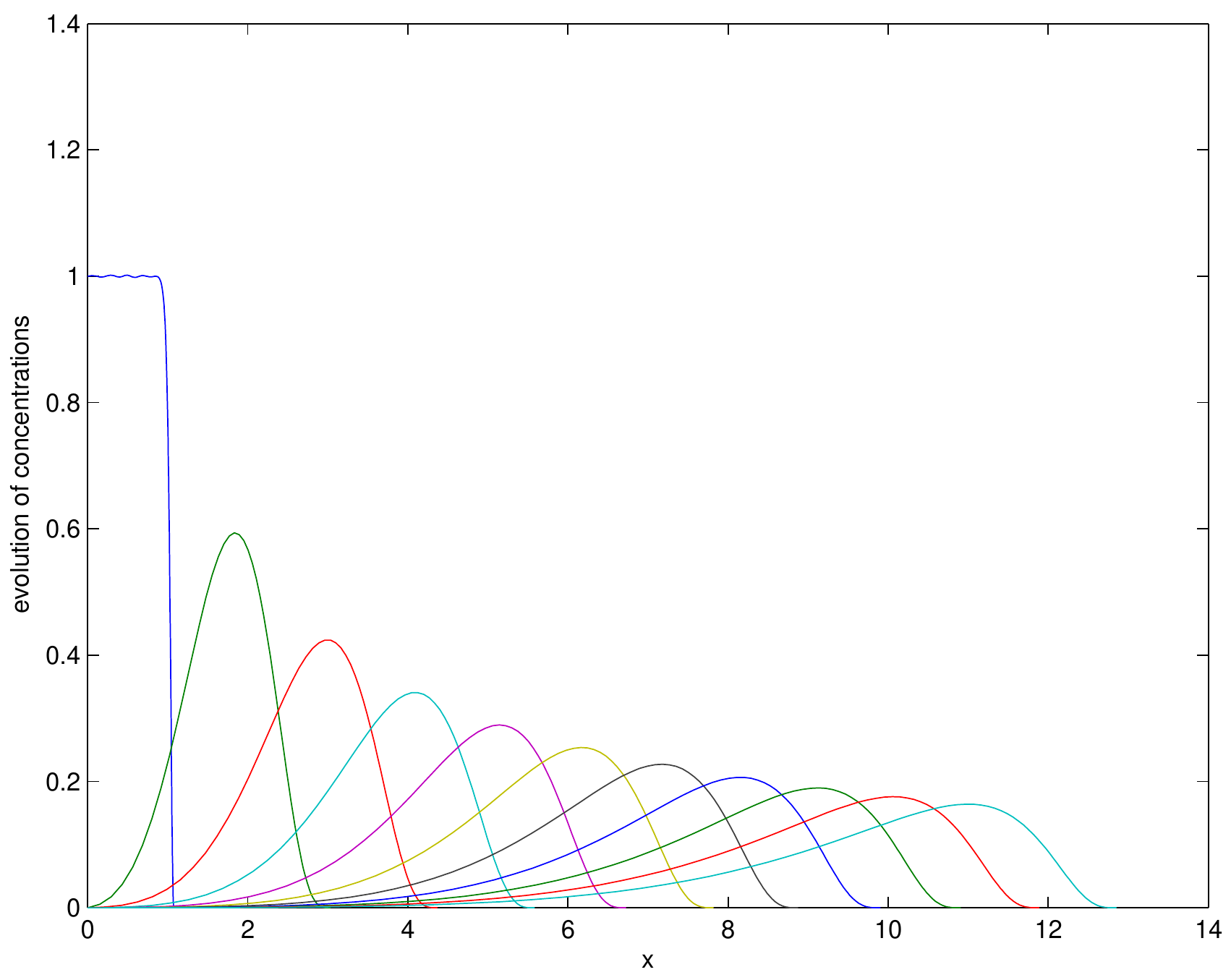}
     }}
 \end{picture}
\caption{Time evolution of concentrations $a=2, p=0.75, N=100, D=0.1, N=100, M=20, t\in (0,40)$}    \label{fig21}
\end{center}
\end{figure}
   \begin{figure}
\begin{center}
  \setlength{\unitlength}{1cm}
\begin{picture}(12,8)
\put(0,0){\mbox{
    \includegraphics[width=12cm,height=6cm] {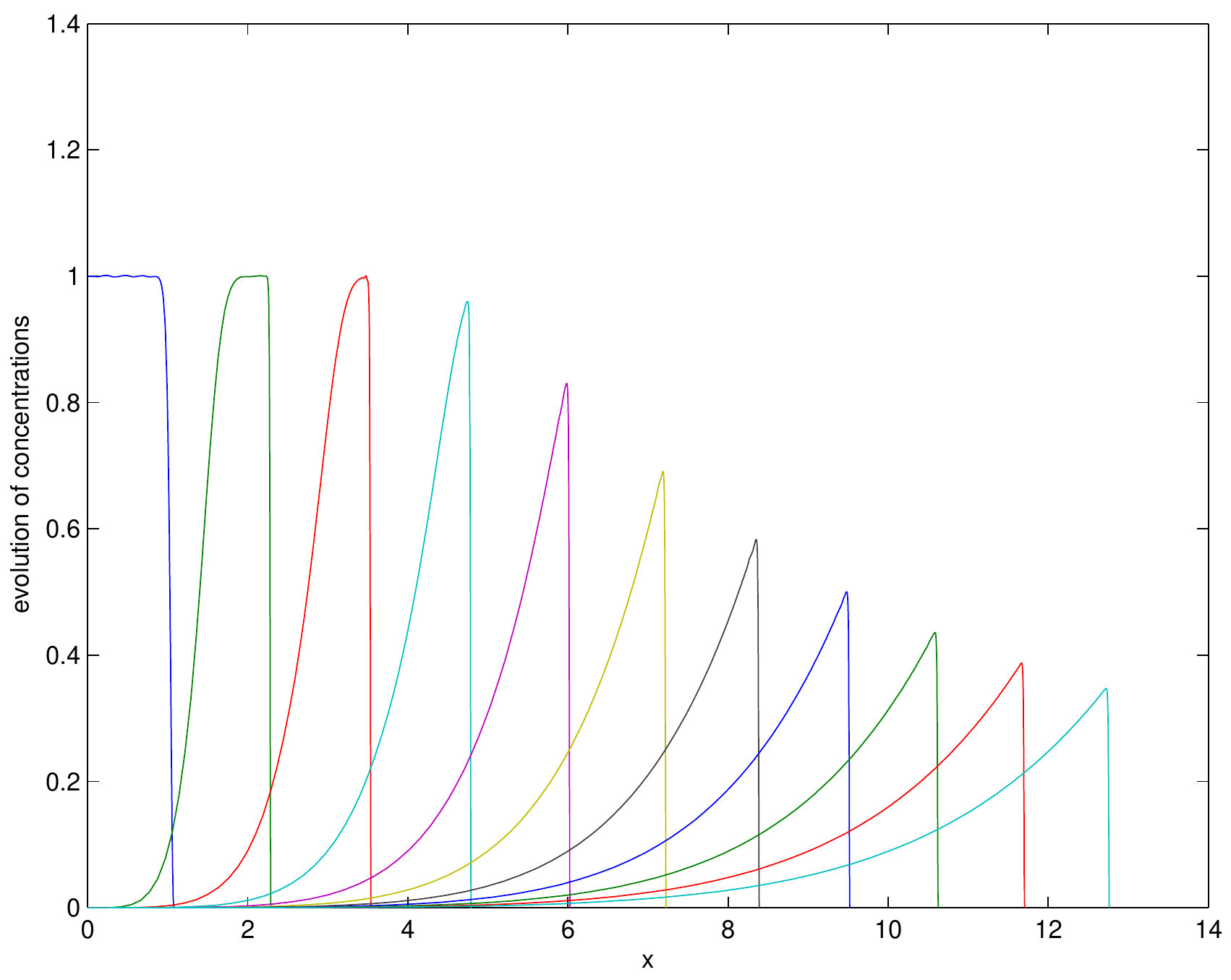}
     }}
 \end{picture}
\caption{Time evolution of concentrations for $a=2, p=0.75, N=100, D=0.001, N=100, M=20, t\in (0,40)$}    \label{fig22}
\end{center}
\end{figure}

\subsection{Simplified model for turbulent flow in porous media}

Consider the following mathematical model
$$ \partial_tu=\partial_x^2u^{3/2}-(u^{3/2}-u^{1/2}) $$
which is a simplified 1D mathematical model for a turbulent flow in porous media -see \cite{[PH]}. Also here
for a special initial profile $u_0$, there is an analytical solution  given by \cite{[PH]}
$$ u(x,t)=\sqrt{a(t)^2+1}[(1-\frac{cosh(x/3)}{\sqrt{a(t)^{-2}+1}})_+]^2 \quad \mbox{with}\quad a(t)=\frac{2(1+\sqrt{2})e^{-5t/6}}{(1+\sqrt{2})^2-e^{-5t/3}}$$
where $(A)_+ =A$ for $A>0$ otherwise $(A)_+=0$. In this case the interface is given by the formula
$$s(t)=3\ln(\sqrt{a(t)^{-2}+1}+a^{-1}(t)).$$
The analytical solution satisfies our interface model given by Theorem 1 (similarly as in Section 7.2).

Following our approximation method in Sections 2-4, we successively obtain $\beta =2$ and the ODE-model for the interface evolution ($x\nearrow s(t)$)
$$ \dot{s}(t) =-3 \partial_xw-\frac {1}{ 2 \partial_x w}, \quad \mbox{where}\ u=w^2.$$
We present our results with discretization parameters $N=60$ and $ M=20$ in the discretization strategy $D_4$. The time evolution of the numerical solution is in Fig. 23. The analytical solution graphically coincides with the  numerical one in the used scaling.
The interfaces are presented in Fig. 24 and the errors in Fig. 25. 
 \begin{figure}
\begin{center}
  \setlength{\unitlength}{1cm}
\begin{picture}(12,8)
\put(0,0){\mbox{
    \includegraphics[width=12cm,height=6cm] {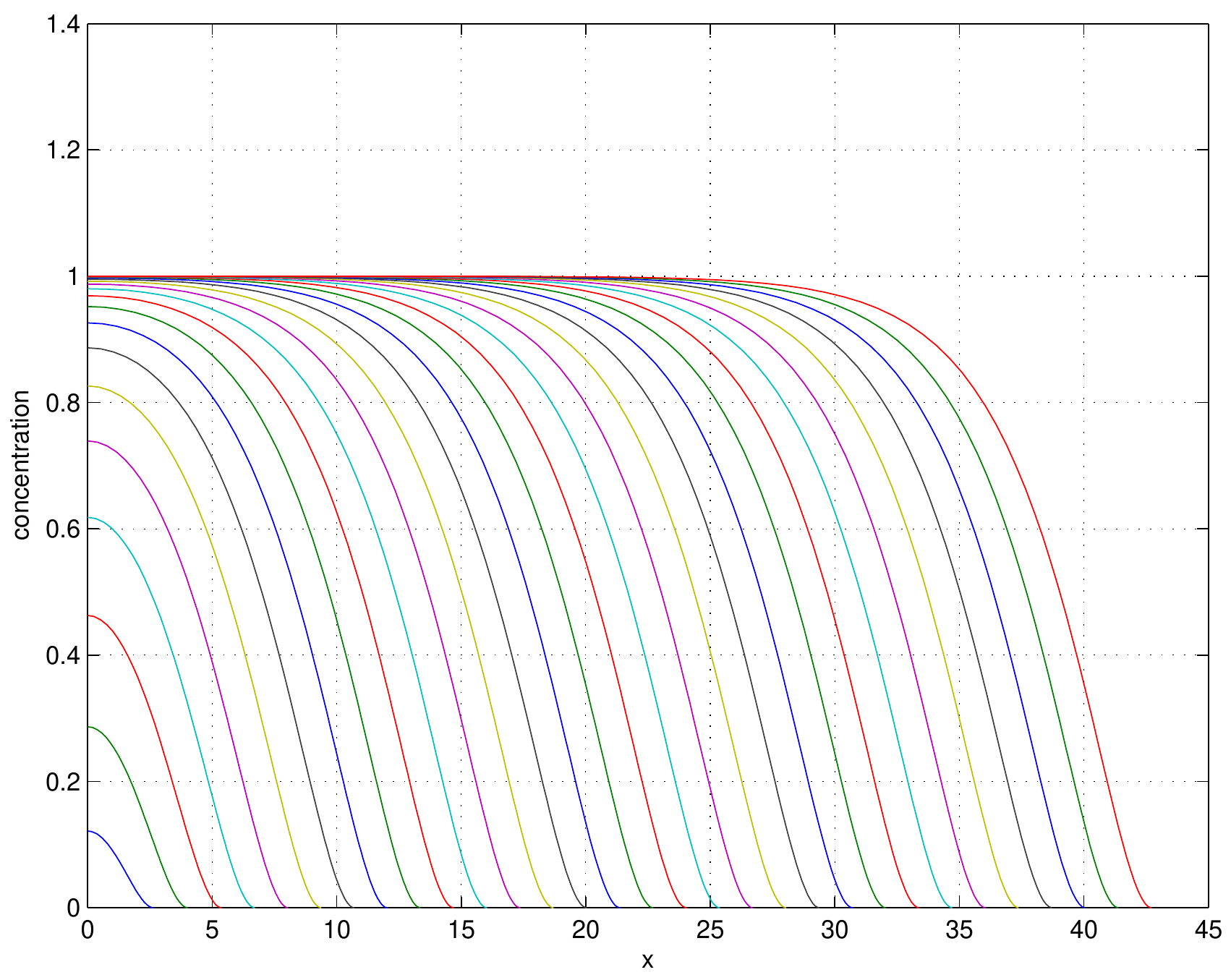}
     }}
 \end{picture}
\caption{Time evolution of the numerical solution for turbulent flow}    \label{fig23}
\end{center}
\end{figure}

 \begin{figure}
\begin{center}
  \setlength{\unitlength}{1cm}
\begin{picture}(12,8)
\put(0,0){\mbox{
    \includegraphics[width=12cm,height=6cm] {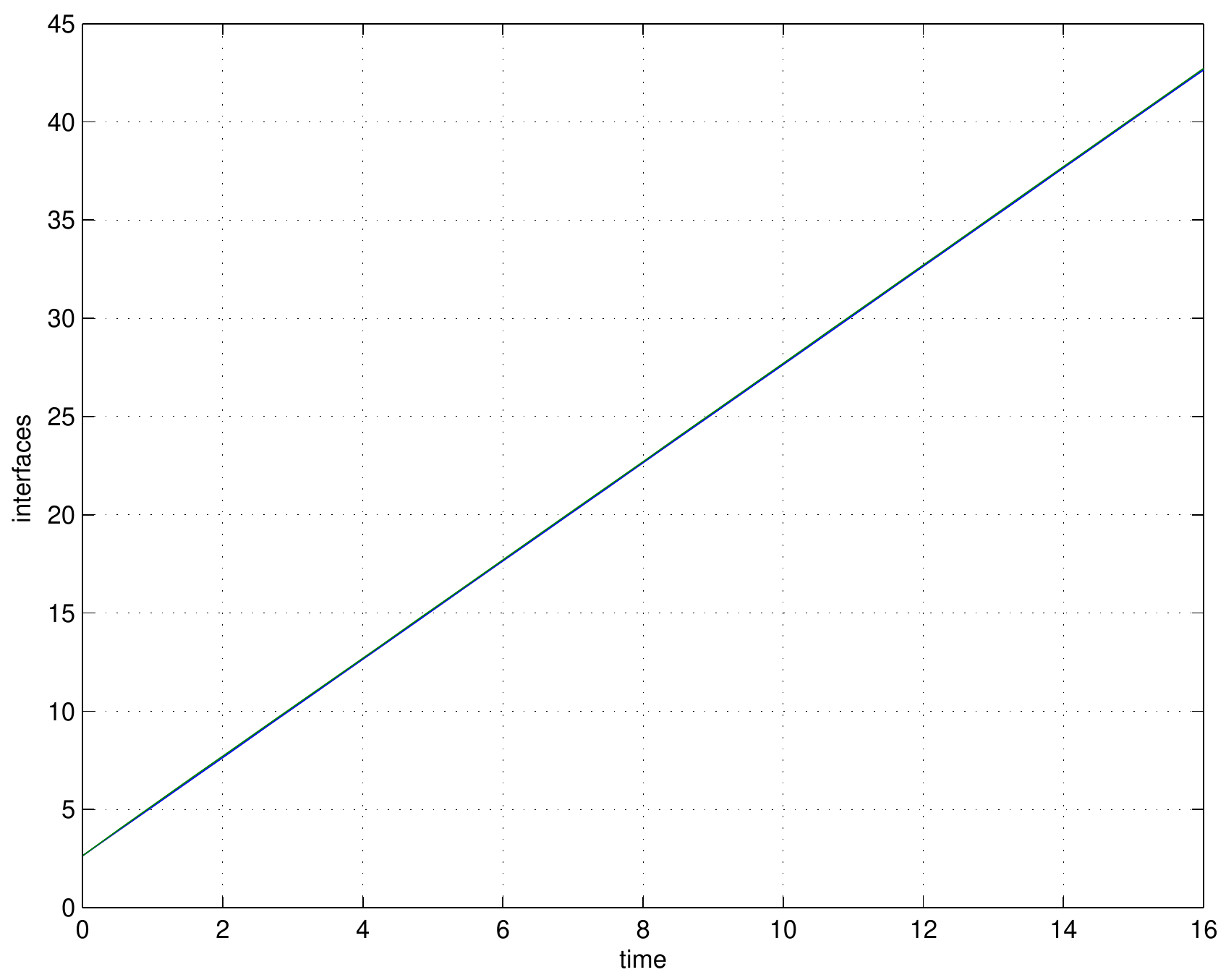}
     }}
 \end{picture}
\caption{Time evolution of interfaces}    \label{fig24}
\end{center}
\end{figure}
 \begin{figure}
\begin{center}
  \setlength{\unitlength}{1cm}
\begin{picture}(12,8)
\put(0,0){\mbox{
    \includegraphics[width=12cm,height=6cm] {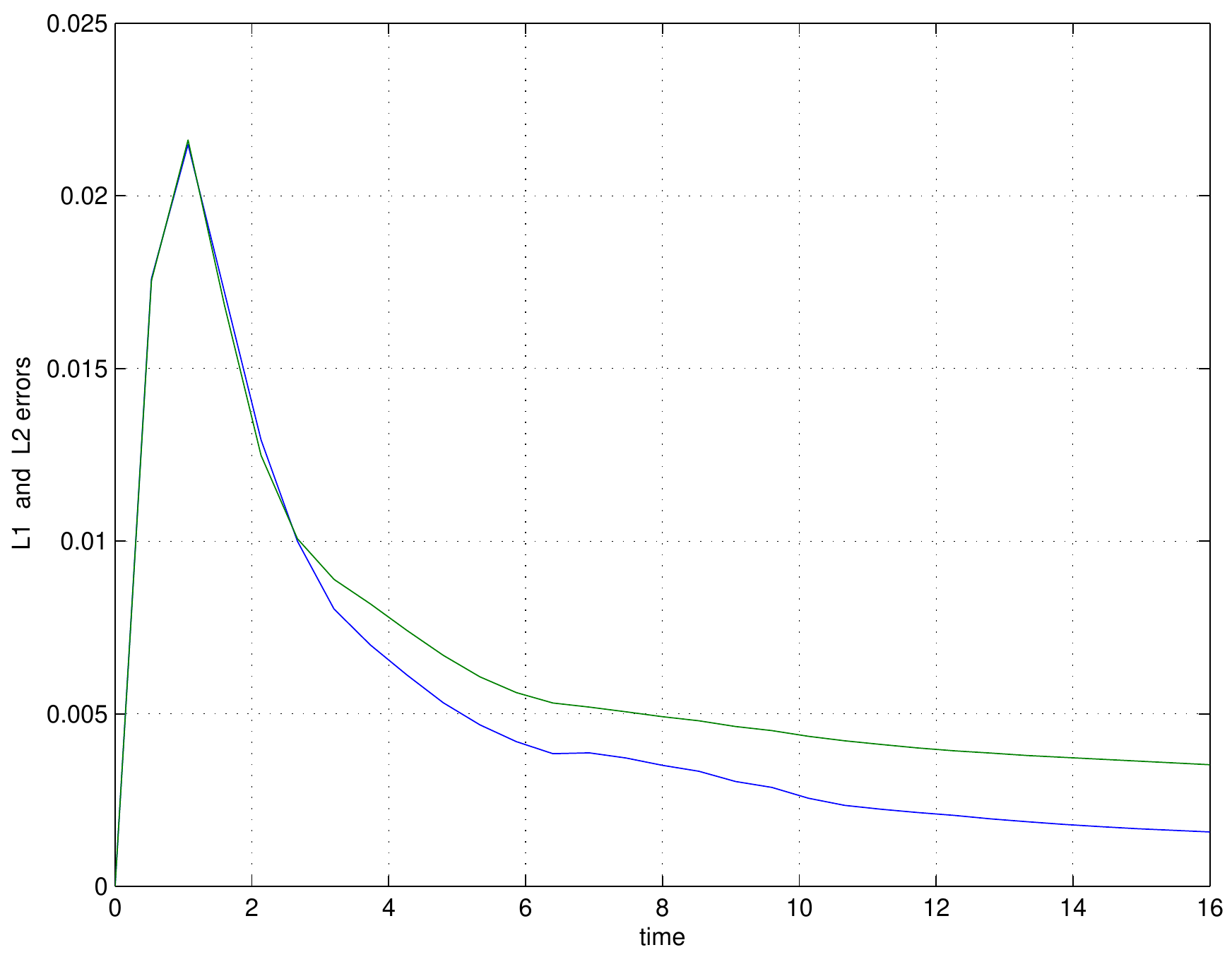}
     }}
 \end{picture}
\caption{Time evolution of errors}    \label{fig25}
\end{center}
\end{figure}
As we can see the numerical and analytical interfaces cannot be distinguished  in the used scaling.
The solution is growing since the source term is positive (the term $u^{1/2}$ dominates $u^{3/2}$).\\
We also present  {\bf the solution of a foam drainage model} 
 (see \cite{[G-K-S],[V-W],[V-W-K]})
 $$ \partial_t u=\partial_x^2 u^{3/2}+\partial_x(u^2) ,\quad \ n>1 .$$
 We consider the boundary condition $\partial_xu|_{x=0}=0$ (symmetry at $x=0$) and the initial condition $u(x,0)=u_0(x)$.\\
 The numerical solution is presented in Fig. 26.
 \begin{figure}
\begin{center}
  \setlength{\unitlength}{1cm}
\begin{picture}(12,8)
\put(0,0){\mbox{
    \includegraphics[width=12cm,height=6cm] {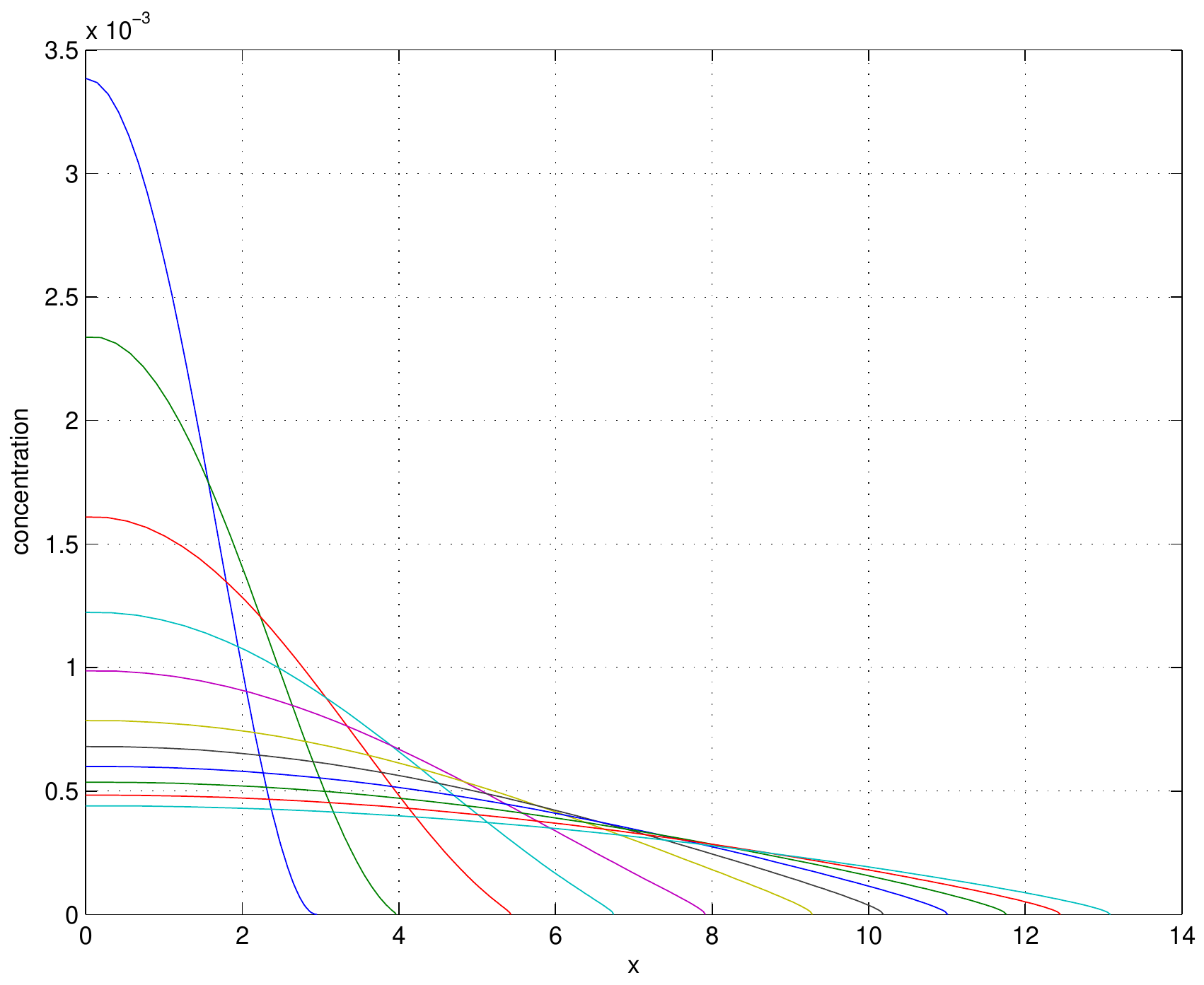}
     }}
 \end{picture}
\caption{Solution of foam drainage model}    \label{fig26}
\end{center}
\end{figure}

Finally, we present the numerical solution for
the mathematical model of viscous liquid advancing over a dry bed.
The motion of a thin sheet of viscous liquid over an inclined plate (see \cite{[Bu]}) is modelled by
 $$ \partial_t u=\partial_x^2 u^4 +\partial_x(u^3),\quad \ n>1 ;$$
 We consider the boundary condition $\partial_xu|_{x=0}=0$ (symmetry at $x=0$) and the initial condition $u(x,0)=u_0(x)$.\\
 The numerical solution is presented in  Fig. 27.
 \begin{figure}
\begin{center}
  \setlength{\unitlength}{1cm}
\begin{picture}(12,8)
\put(0,0){\mbox{
    \includegraphics[width=12cm,height=6cm] {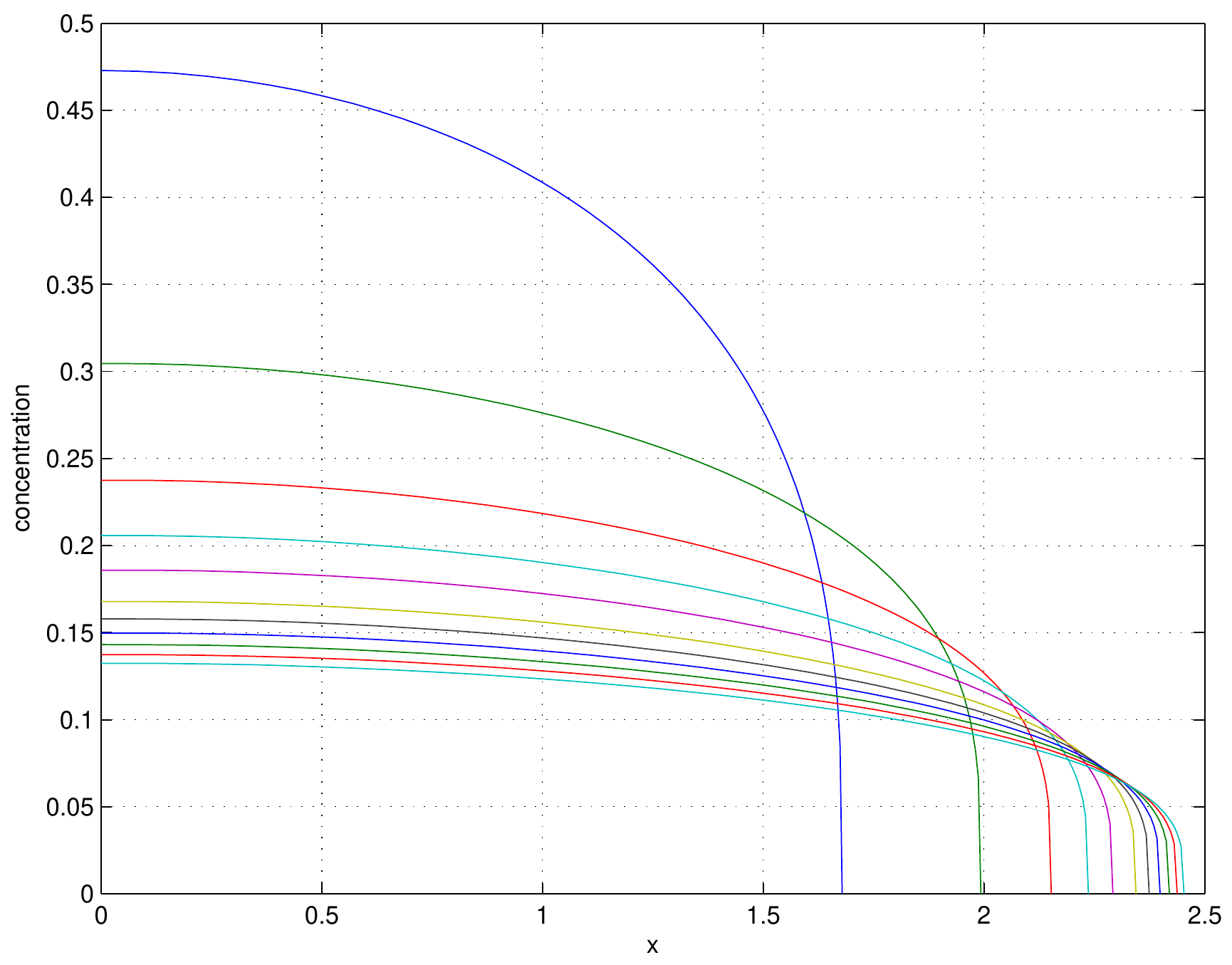}
     }}
 \end{picture}
\caption{Advancing of viscous liquid}    \label{fig27}
\end{center}
\end{figure}

\bibliographystyle{elsart-num-sort}

\end{document}